\documentclass[11pt]{elsarticle}
\usepackage{hyperref}
\usepackage{pdfpages}

\journal{Journal of Econometrics}

\usepackage{color, amsmath, amssymb, amsthm, graphicx, hyperref}
\usepackage[top=1in, bottom=1in, left=1in, right=1in]{geometry}
\usepackage[toc,page]{appendix}
\hypersetup{colorlinks=true, citecolor=blue}

\allowdisplaybreaks

\usepackage{setspace}
\usepackage{xr}
\usepackage{etoolbox}

\newtheorem{theorem}{Theorem}
\newtheorem{definition}{Definition}
\newtheorem{proposition}{Proposition}
\newtheorem{corollary}{Corollary}
\newtheorem{lemma}{Lemma}
\newtheorem{remark}{Remark}
\newcommand{\possessivecite}[1]{\citeauthor{#1}'s (\citeyear{#1})}
\newcommand{\possessivecitewos}[1]{\citeauthor{#1}' (\citeyear{#1})}

\biboptions{authoryear}

\begin{document}

\begin{frontmatter}

\title{On the Predictive Risk in Misspecified Quantile Regression}

\author[mymainadress]{Alexander Giessing}
\ead{giessing@umich.edu}
\author[mymainadress]{Xuming He\corref{mycorrespondingauthor}}
\ead{xmhe@umich.edu}

\address[mymainadress]{Department of Statistics, University of Michigan, 1085 South University Avenue, Ann Arbor, Michigan 48109, U.S.A. }

\cortext[mycorrespondingauthor]{Corresponding author}

\begin{abstract}
	In the present paper we investigate the predictive risk of possibly misspecified quantile regression functions. The in-sample risk is well-known to be an overly optimistic estimate of the predictive risk and we provide two relatively simple (asymptotic) characterizations of the associated bias, also called expected optimism. We propose estimates for the expected optimism and the predictive risk, and establish their uniform consistency under mild conditions. Our results hold for models of moderately growing size and allow the quantile function to be incorrectly specified. Empirical evidence from our estimates is encouraging as it compares favorably with cross-validation.
\end{abstract}

\begin{keyword}
Quantile Regression \sep Misspecification \sep Predictive Risk \sep Expected Optimism
\JEL C14 \sep C51 \sep C52 \sep C53
\end{keyword}

\end{frontmatter}

\section{Introduction}\label{sec:Intro}
Predictive modeling is at the core of many scientific disciplines, including business, engineering, finance, and public health. A natural way to gauge the predictive capability of a statistical model is to estimate its predictive risk. The systematic study of the risk of a statistical procedure traces back to at least~\cite{stein1981}. Since then, the concept of risk has become an integral part of applied statistical modeling: predictive risk is routinely used to assess the complexity of statistical modeling procedures~\citep[e.g.][]{akaike1973, mallows1973, foster1994} to compare different statistical models~\citep[e.g.][]{hastie1990, ye1998}, and to choose between tuning parameters that control bias-variance trade-offs~\citep[e.g.][]{donoho1995, kou2002}. In several special cases,~\possessivecite{stein1981} theory of unbiased risk estimation provides simple estimates for the risk of a statistical model. However, in general, there does not exist a unified approach to estimating the predictive risk of a statistical model or procedure. 

In this paper, we focus on the predictive risk of possibly misspecified quantile regression models. In addition to its role in applied statistical modeling as outlined above, the predictive risk from quantile regression models has also garnered significant interest in finance and risk management to assess the value-at-risk and expected shortfall of investments or to solve portfolio choice problems~\citep[e.g.][]{xiao2015,cahuich2013, gaglianone2011,hexd2011,bassett2004, engle2004, chernozhukov2001}.

We contribute to the theory on the predictive risk of quantile regression models by deriving two asymptotic characterizations of the bias of the in-sample risk (when used to estimate the predictive risk) and proposing a uniformly consistent, de-biased estimator for the predictive risk. Following the terminology introduced by~\cite{efron1983} we call the bias of the in-sample risk the ``expected optimism''.

Our first characterization of the expected optimism provides a characterization comparable to~\possessivecite{efron2004} covariance penalty and~\possessivecite{tibshirani1999} covariance inflation criterion. The second characterization relates to robust and generalized Akaike-type information criteria~\citep[e.g.][]{lv2014, portnoy1997, burman1995}. Both characterizations show that large part of the expected optimism can be attributed to a nonlinear function of the quantile level, the conditional density of the response variable given the predictors and the (weighted) covariance matrix of the predictors. Specializing to location models, we glean additional insight into the expected optimism and its functional dependence on the conditional density and the number of predictors. As a consequence, the commonly used notion of effective degree of freedom for a statistical model has a richer content for misspecified models.

The second characterization of the expected optimism lends itself to a simple plug-in estimator. We establish its uniform consistency over a class of candidate models and, based on this result, propose a uniformly consistent, de-biased estimator of the predictive risk. Our theoretical analysis indicates that the de-biased estimator is particularly relevant in the case in which the dimension of candidate models grows at least in the order of the square root of the sample size. Empirical evidence suggests that the de-biasing procedure is practically relevant even when the model size is fixed and relatively small compared to the sample size. A comparison of our de-biased estimate against the popular method of cross-validation is favorable for our procedure.

To allow broad applicability our our results, we develop our theory in a triangular array of row-wise independent random vectors whose dimension may grow with the sample size. We only require minimal assumptions on the joint distribution of the response and predictor variables. Notably, the response and the predictor variables can both be unbounded, their marginal distributions can be non-Gaussian, and their relationship (i.e. the conditional quantile functions) can be linear, nonlinear or nonparametric. Thus, our framework for quantile regression generalizes the frameworks of~\cite{lee2015, noh2013, angrist2006, kim2003} who consider misspecified quantile regression models with a fixed number of parameters. Unlike the recent literature on quantile regression based on series, semi- and nonparametric estimators we do not assume that the misspecification error vanishes as more predictors are included in the regression function~\citep{belloni2016, chao2016}. Naturally, our results continue to hold if the model is (asymptotically) correctly specified.

We organize this article as follows: In Section~\ref{sec:MisspecifiedQRandPredRisk} we lay out a general framework for misspecified quantile regression models. We introduce necessary terminology and discuss how to define the predictive risk of potentially misspecified quantile regression models. In Section~\ref{sec:ApproxExpectOpt} we derive two asymptotic characterizations of the expected optimism of the in-sample risk and discuss insights that we gain from these characterizations. In Section~\ref{sec:EstimationPredRisk} we propose a nonparametric plug-in estimator for one of the asymptotic characterizations of the expected optimism and use it to construct a de-biased estimate of the predictive risk. We establish uniform consistency of both estimators. In Section~\ref{sec:NumericalEvidence} we report numerical evidence that our estimates of the expected optimism and the predictive risk are on target, and that the predictive risk estimate can be better than the commonly-used cross-validation approach. We conclude in Section~\ref{sec:Conclusion} with additional remarks, and present all proofs in~\ref{appendix:Proofs}. We provide further theoretical results and technical lemmata in the Supplementary Materials~\ref{Supplement:Angrist},~\ref{Supplement:ModelSelection}, and~\ref{Supplement:TechnicalLemmata}. Additional simulation results are relegated to the Supplementary Materials~\ref{Supplement:NumericalEvidence}.

\section{Misspecified quantile regression and predictive risk}\label{sec:MisspecifiedQRandPredRisk}
\subsection{Notation and framework}\label{subsec:NotationFramework}
The setting of interest is a high-dimensional triangular array $\mathcal{D}_n = \{(Y_{ni},X_{ni})\}_{i=1}^n$, where $(Y_{ni},X_{ni}) \in \mathbb{R} \times \mathcal{X}$ are row-wise independent random vectors with distribution $F_n$ which may change with the sample size $n$. As per convention the scalar variable $Y_{ni}$ denotes the response variable and the vector $X_{ni} \in \mathcal{X}$ denotes a vector of covariates. We denote by $F_{Y_n|X_n}$ the conditional distribution of $Y_{ni}$ given $X_{ni}$. We use subscripts on the expectation operator $\mathbb{E}$ to specify to which random variable the operator is applied to, i.e. $\mathbb{E}_{(Y_{n1},X_{n1})}$ means that expectation is only taken over $(Y_{n1},X_{n1})$ whereas $\mathbb{E}_{\mathcal{D}_n}$ means that expectation is taken over the entire triangular array $\mathcal{D}_n$. We let
\begin{align}\label{eq:XZMap}
x \mapsto Z(x) = \big(Z_1(x), \ldots, Z_d(x)\big)
\end{align}
denote a mapping from $\mathcal{X}$ into $\mathbb{R}^d$ and call the transformed covariates $Z(X_{n1}), \ldots, Z(X_{nn})$ predictor variables. We consider the case where the dimension $d$ of the predictor variables grows with the sample size $n$ and may be much larger than $n$. We call a subset $S \subseteq \{1, \ldots, d\}$ of predictors $Z(X_{ni})$ a model and write
\begin{align}
Z_S(X_{ni}) = \big(Z_j(X_{ni})\big)_{j \in S}.
\end{align}
We denote the collection of models under consideration by $M$. We allow $M$ to be as large as the power set of $\{1, \ldots, d\}$ and to grow with the sample size $n$. We write $|S|$ for the cardinality of a model $S$ and denote the largest cardinality of models in $M$ by $m$. Clearly, we have $m \leq d$. 

The purpose of linear quantile regression is to approximate the true conditional quantile function (CQF) of $Y_{ni}$ given $X_{ni}$, 
\begin{align}\label{CQF}
Q_{Y_n}(\tau|X_{ni}) = \inf\left\{y: F_{Y_{n}|X_{n}}(y|X_{ni}) \geq \tau \right\},
\end{align}
by a linear function of the predictor variables $Z(X_{ni})$. To this end, we assume that the vectors of predictor variables $Z(X_{ni})$ consist of series functions with reasonably good approximation properties such as indicators, B-splines, regression splines, polynomials, Fourier series, or wavelets~\citep[e.g.][]{belloni2016, chao2016}. However, unlike them we do not require that the approximation error vanishes as the number of predictors $m$ increases, i.e. we allow for persistent misspecification. We define the vector of regression coefficients $\theta^\tau_{n,S} = (\theta_{n1}^\tau, \ldots, \theta_{n|S|}^\tau)'$ associated with model $S$ as the solution to the quantile regression problem 
\begin{align}\label{eq:PopulationQRProblem}
\min_{\theta \in \mathbb{R}^{|S|}} \mathbb{E}_{\mathcal{D}_n} \left[\rho_\tau\big(Y_{n1} - Z_S(X_{n1})'\theta\big) - \rho_\tau\big(Y_n - Q_{Y_n}(\tau|X_{n1})\big) \right],
\end{align}
and the vector of estimated regression coefficients $\hat{\theta}_{n,S}^\tau = (\hat{\theta}_{n1}^\tau, \ldots, \hat{\theta}_{n|S|}^\tau)'$ as the solution to the sample quantile regression problem
\begin{align}\label{eq:SampleQRProblem}
\min_{\theta \in \mathbb{R}^{|S|}} \frac{1}{n} \sum_{i=1}^n \rho_\tau\big(Y_{ni} - Z_S(X_{ni})'\theta\big),
\end{align}
where $\rho_\tau(u) = (\tau - 1\{u \leq 0\})$ is the check loss~\citep{koenker2005}. The estimate of the true CQF of $Y_{n}$ given $X_{n}$ based on model $S$ is given as
\begin{align}
\widehat{Q}_{Y_n}(\tau|X_{n}, S) = Z_S(X_{n})'\hat{\theta}_{n,S}^\tau.
\end{align}

\cite{koenker1978} show that under mild conditions $\widehat{Q}_{Y_n}(\tau|X_n,S)$ is a consistent estimate of $Q_{Y_n}(\tau|X_n)$ if the true CQF is indeed linear in $Z_S(X_{ni})$ and the dimension of the predictor variables is fixed.~\cite{angrist2006} establish corresponding consistency and asymptotic normality results for misspecified quantile regression models. The results on general M-estimators~\citep[][]{he2000} (Theorem 1), semi-parametric quantile regression~\citep[][]{chao2016} and quantile series estimators~\citep[][]{belloni2016} extend these result to cases in which the dimension of the predictors $m$ increases with the sample size $n$.

\subsection{Predictive risk and expected optimism}\label{subsec:PredRiskandExpectOptGENERAL}
Two statistical theories have been developed to estimate the predictive risk, cross-validation~\citep[e.g.][]{stone1974, stone1977, allen1974, golub1979, wahba1990, efron1983, efron1986, efron2004, efron1997} and covariance penalties, which  include techniques such as~\possessivecite{mallows1973} Cp, \possessivecite{akaike1969} information criterion (AIC) and final prediction error (FPE), \possessivecite{takeuchi1976} information criterion (TIC), and~\possessivecite{stein1981} unbiased risk estimate (SURE). Our approach to estimating the predictive risk of potentially misspecified quantile regression models falls into the category of covariance penalties. In this section we therefore introduce necessary terminology and the rational behind covariance penalties.

Suppose that a model $f$ is fitted to some data $\mathcal{Z}_n = \{Z_1, \ldots, Z_n\}$ producing an estimate $\hat{\mu}_n = f(\mathcal{Z}_n)$ for target $\mu$. Predictive risk evaluation tries to assess how well $\hat{\mu}_n$ predicts $\mu$ at a future data point $Z^0$ independently generated from the same mechanism that produced $\mathcal{Z}_n$. To measure the error between $\hat{\mu}_n$ and $\mu$ one chooses a loss function $L$ and defines the predictive risk as the average loss over current and future data, i.e.
\begin{align}
\mathbb{E}_{\mathcal{Z}_n, Z^0}\left[L\Big(\mu(Z^0), \hat{\mu}_n(Z^0)\Big)\right].
\end{align}

Covariance penalties provide as an intermediate result an estimate of the bias of the in-sample risk when used as estimate of the predictive risk. Following the terminology introduced by~\cite{efron1983} we call the negative bias the ``expected optimism'' of the in-sample risk, 
\begin{align}
b_n(L, \mu) = \mathbb{E}_{\mathcal{Z}_n, Z^0}\left[L\Big(\mu(Z^0), \hat{\mu}_n(Z^0)\Big)\right] - \mathbb{E}_{\mathcal{Z}_n}\left[\frac{1}{n}\sum_{i=1}^nL\Big(\mu(Z_i), \hat{\mu}_n(Z_i)\Big)\right].
\end{align}

Given a consistent estimate $\hat{b}_n(L, \mu)$ of $b_n(L, \mu)$ one obtains a consistent and de-biased estimate of the predictive risk via
\begin{align}
\frac{1}{n}\sum_{i=1}^nL\Big(\mu(Z_i), \hat{\mu}_n(Z_i)\Big) + \hat{b}_n(L, \mu).
\end{align}

Even though covariance penalties are conceptually straightforward, so far, they have only been derived for a limited number of loss functions, namely the square loss and the ``q class of error functions''~\citep[][]{efron1986}. This is likely because the expected optimism of most other loss functions is a highly non-linear function for which it is difficult to construct estimators. Yet, in principle, covariance penalties have two advantages over cross-validation techniques: First, cross-validation techniques tend to produce estimates of the predictive risk that have a high variance than covariance penalties, since they split the sample into test and training sets and thereby reduce the number of samples from which $\hat{\mu}_n$ is estimated~\citep[e.g.][]{efron2004}. Second, cross-validation techniques are known to produce biased estimates of the predictive risk. Several heuristic adjustments to (the vanilla) cross-validation techniques have been proposed, but they lack rigorous proofs~\citep[e.g.][]{burman1989, tibshirani2009}.

\subsection{Predictive risk and expected optimism in quantile regression}\label{subsec:PredRiskandExpectOptQR}
We discuss the choice of the loss function to measure the predictive risk of a potentially misspecified quantile regression model $S$ and define the associated expected optimism.

Let $(Y_n^0, X_n^0)$ be a pair of data points drawn from $F_n$ and independent of sample $\mathcal{D}_n=\{(Y_{ni}, X_{ni})\}_{i=1}^n$. Fix a model $S \subseteq \{1, \ldots, d\}$ and consider the estimate of the CQF of $Y_n^0$ given $X_n^0$ based on model $S$ and sample $\mathcal{D}_n$, i.e.
\begin{align}
\widehat{Q}_{Y_n^0}(\tau|X_n^0, S) = Z_S(X_n^0)'\hat{\theta}_{n,S}^\tau.
\end{align}

Since the true CQF of $Y^0_n$ given $X^0_n$, $Q_{Y_n^0}(\tau|X_n^0)$, is not an observable statistic given the data $\mathcal{D}_n$ and $(Y^0_n, X^0_n)$, risk measures which assess directly the difference between estimate $\widehat{Q}_{Y_n^0}(\tau|X_n^0, S)$ and target $Q_{Y_n^0}(\tau|X_n^0)$, such as the mean squared prediction error or the mean absolute prediction error, do not have (simple) sample analogues. We therefore propose the following risk measure which depends only on observables.
\begin{definition}[Predictive risk]\label{def:PredRisk} The predictive risk of quantile regression model $S$ is 
	\begin{align*}
	\mathrm{PR}^\tau_n(S) &= \mathbb{E}_{\mathcal{D}_n, (Y_n^0, X_n^0)}\left[\rho_\tau\big(Y^0_n - \widehat{Q}_{Y_n^0}(\tau|X_n^0, S)\big) - \rho_\tau(Y^0_n)\right],
	\end{align*}
	where $(Y_n^0, X_n^0)$ is a pair of data points drawn from $F_n$ and independent of sample $\mathcal{D}_n$.
\end{definition}

The associated expected optimism of using the in-sample risk $\frac{1}{n} \sum_{i=1}^n \Big(\rho_\tau\big(Y_{ni} - \widehat{Q}_{Y_{n}}(\tau|X_{ni}, S)- \rho_\tau(Y_{ni})\Big)$ as an estimate of the predictive risk is defined as follows.
\begin{definition}[Expected Optimism]\label{def:ExpectedOptimism} The expected optimism of quantile regression model $S$ is 
\begin{align*}
b_n^\tau(S) &= \mathrm{PR}^\tau_n(S) - \mathbb{E}_{\mathcal{D}_n}\left[\frac{1}{n} \sum_{i=1}^n \Big(\rho_\tau\big(Y_{ni} - \widehat{Q}_{Y_{n}}(\tau|X_{ni}, S)\big) - \rho_\tau(Y_{ni})\Big)\right].
\end{align*}
\end{definition}

Several comments are in order with regard to these two definitions. First, the reason for subtracting $\rho_\tau(Y^0_n)$ in Definition~\ref{def:PredRisk} (and $\rho_\tau(Y_{ni})$ in Definition~\ref{def:ExpectedOptimism}) is purely technical: it allows us to dispense with moment conditions on the response variable $Y^0_n$. To see this, note that the check loss $\rho_\tau$ is Lipschitz continuous and hence the predictive risk $\mathrm{PR}^\tau_n(S)$ is upper bounded by $\mathbb{E}_{\mathcal{D}_n, (Y_n^0, X_n^0)}\big|Z_S(X^0_n)'\hat{\theta}^\tau_{n,S} \big|$. For this expected value to be finite it suffices that the CQF of $Y^0_n$ given $X^0_n$ has finite second moments~\citep[e.g.][]{angrist2006}.

Second, the predictive risk of model $S$ can be shown to be an (almost) affine transformation of the unconditional mean squared prediction error (MSPE) $\mathbb{E}_{\mathcal{D}_n, X_n^0} \big[\big(Q_{Y_n^0}(\tau|X_n^0) - \widehat{Q}_{Y_n^0}(\tau|X_n^0, S)\big)^2\big]$, which cannot be estimated directly since it depends on $Q_{Y_n^0}(\tau|X_n^0)$, the unobserved true CQF. Since the MSPE is itself an important quantity to assess model fit, its connection with our notion of predictive risk may be of independent interest. We relegate the precise statement of this technical result to the Supplementary Materials~\ref{Supplement:Angrist}.

Third, the predictive risk based on the check loss $\rho_\tau$ has garnered significant interest in finance and risk management. For example, it is used in the context of value-at-risk~\citep[e.g.][]{xiao2015, gaglianone2011}, conditional value-at-risk and expected shortfall~\citep[e.g.][]{engle2004, chernozhukov2001} and portfolio choice problems with Choquet expectation~\cite[e.g.][]{cahuich2013, hexd2011, bassett2004, tversky1992}.

Fourth, the predictive risk and the expected optimism play an important role in model selection criteria~\citep[e.g.][]{akaike1973, ronchetti1985, foster1994, burman1995, portnoy1997, ye1998, bozdogan2000, lv2014}, model comparison~\citep[e.g.][]{hastie1990, tibshirani1999, kou2002}, and computation of generalized degrees of freedom~\citep[e.g.][]{ye1998}.

\subsection{Technical assumptions}\label{sec:Assumptions}
For the theoretical investigations of the predictive risk and the expected optimism of potentially misspecified quantile regression models we require several assumptions, which we discuss in this section. Since the quantile level $\tau$ is always pre-specified, we suppress the dependence on $\tau$ in some notation. Recall that $S \subseteq \{1, \ldots, d\}$, $|S| \leq m$, and that $M$ is a subset of the power set of $\{1, \ldots, d\}$. Throughout, we assume that $M$ contains at least two models, i.e. $|M| \geq 2$, and that $n \geq 16$, i.e. $\log \log n \geq 1$. 

\textit{
\begin{itemize}
	\item[(A1)] The data $(Y_{ni}, X_{ni}) \in \mathbb{R} \times \mathcal{X}$ are row-wise independent random vectors with distribution $F_n$, where $F_n$ may change with the sample size $n$.
	\item[(A2)] The conditional density $f_{Y_n|X_n}$ of $Y_n$ given $X_{n}$ is uniformly bounded from above, i.e. there exits $\nu_{+} < \infty$ such that
	\begin{align*}
	\limsup_{n \rightarrow \infty}\sup_{a \in \mathbb{R}}\sup_{x \in \mathbb{R}^d}\left| f_{Y_n|X_n}(a|x)\right| \leq \nu_{+}.
	\end{align*}
	\item[(A3)] The conditional density $f_{Y_n|X_n}$ of $Y_n$ given $X_n$ is $\alpha$-H{\"o}lder continuous for $\alpha \in \left[\frac{1}{2}, 1\right]$, i.e. there exists a constant $\nu_H > 0$ such that for any $a, b \in \mathbb{R}$,
	\begin{align*}
	\limsup_{n \rightarrow \infty}\sup_{x \in \mathbb{R}^d}\Big| f_{Y_n|X_n}(a|x) - f_{Y_n|X_n}(b|x)\Big| \leq \nu_{H} |a - b|^{\alpha}.
	\end{align*}
	\item[(A4)] The maximum eigenvalue of the matrix of second moments is uniformly bounded from above, i.e. there exists $\lambda_{+} < \infty$ such that
	\begin{align*}
	\limsup_{n \rightarrow \infty} \max_{S \in M}\lambda_{\max } \left(\mathbb{E}_{X_n}\left[Z_S(X_{n})Z_S(X_{n})'\right] \right) \leq \lambda_+,
	\end{align*} and  the minimum eigenvalue of the weighted second moment matrix is bounded from below by $\lambda_{n} > 0$,
	\begin{align*}
	\min_{S \in M}\lambda_{\min}\left(\mathbb{E}_{X_n}\left[f_{Y_n|X_n}\big(Z_S(X_{n})'\theta^\tau_{n,S}|Z_S(X_n)\big)Z_S(X_{n})Z_S(X_{n})'\right]\right) > \lambda_{n}.
	\end{align*}
\end{itemize}
}

In the above assumptions the uniformity in $n$ is necessary since we consider triangular arrays. Assumptions (A1), (A2), and (A3) with $\alpha =1$ are fairly standard in the quantile regression literature~\cite[e.g.][]{angrist2006,belloni2016,chao2016}. It is possible to relax the (implicit) assumption that the random variables are identically distributed within each row; in fact independence suffices for our results. However, we do not pursue these refinements in the present paper. The stringentness of Assumption (A4) depends on how fast $\lambda_n$ is allowed to go to zero. We require the following technical rate condition on $\lambda_n$:

\textit{
\begin{itemize}
	\item[(A5)] The minimum eigenvalue of the matrix of second moments, $\lambda_n$, is bounded below asymptotically in the following way:
	\begin{align*}
	\lambda_n \gtrsim \left(\frac{m \log |M| \: \log \log n}{n}\right)^{1/2 - 1/(4\alpha)}.
	\end{align*}
\end{itemize}
}

This rate condition is purely technical and difficult to motivate. Clearly, the condition is less stringent the larger $\alpha$, i.e. the smoother the conditional density $f_{Y_n|X_n}$ of $Y_n$ given $X_n$. In particular, if $\alpha = 1/2$, we require $\lambda_n = O(1)$; whereas in the case of a continuous conditional density, we allow $\lambda_n = O\big((m \log |M| \: \log \log n)^{1/4}n^{-1/4}\big)$. The rate condition relaxes the stronger boundedness assumptions on the largest and smallest eigenvalue of the weighted second moment matrix that prevail in the literature on quantile regression~\cite[][]{koenker2017}. Together with the upper bound on the largest eigenvalue of the expected value of the Gram-matrix the rate condition implies that $m \lesssim n$. This is a much weak condition on the growth rate of the number of predictors than has been proposed in recent work on (misspecified) quantile regression with increasing number of predictors. E.g.~\cite{belloni2016} and~\cite{chao2016} require that $\zeta_m \equiv \sup_{x \in \mathcal{X}}\|Z(x)\|_2 < \infty$ satisfies $m \zeta_m^2 (\log n)^2 = o(n)$. If the predictors are element-wise bounded, this amounts to the condition $m^2 (\log n)^2 = o(n)$. We shall see that our relaxed assumption on the growth rate is important in the theoretical analysis of the proposed estimate for the predictive risk in Section~\ref{sec:EstimationPredRisk}.

Lastly, we introduce the following moment condition on the predictors:

\textit{
\begin{itemize}
	\item[(A6)] The vector $Z(X_n) = \big(Z_1(X_{n}), \ldots, Z_d(X_{n})\big)$ is a vector of random variables with finite $8 + \delta$ moment, for some $\delta > 0$. In particular, for $1 \leq k \leq 8$, there exist constants $\mu_k > 0$ such that
	\begin{align*}
	\limsup_{n \rightarrow \infty} \max_{j =1, \ldots, d}\left(\mathbb{E}_{X_n}\left[\left|Z_j(X_{n}) \right|^{k + \delta} \right]\right)^{1/(k + \delta)} \leq \mu_k.
	\end{align*}
\end{itemize}
}
This condition is significantly weaker than the uniform boundedness assumption on the map $Z$ imposed in~\cite{belloni2016} and~\cite{chao2016} (i.e. $\zeta_m \equiv \sup_{x \in \mathcal{X}}\|Z(x)\|_2 < \infty$). Again, uniformity in $n$ is necessary since we consider triangular arrays.

\section{Two asymptotic characterizations of the expected optimism}\label{sec:ApproxExpectOpt}
\subsection{The covariance form of the expected optimism}\label{subsec:CovExpectOpt}
In the case of ordinary least squares, the expected optimism can be evaluated via~\possessivecitewos{mallows1973} $C_p$. In the case of nonlinear least squares with Gaussian errors, the expected optimism can be estimated via~\possessivecite{stein1981} divergence formula. And for loss functions that belong to~\possessivecite{efron2004} ``q class of error measures'', the expected optimism can be expressed as a function of the covariance of two observable quantities.

Since the expected optimism $b^\tau_n(S)$ from Definition~\ref{def:ExpectedOptimism} is based on the check loss $\rho_\tau$, none of the above three results applies. Instead we have the following result.

\begin{theorem}[Covariance Form of the Expected Optimism]\label{theorem:ExpectOptCovForm}
	Suppose that Assumptions (A1) -- (A6) from Section~\ref{sec:Assumptions} hold. Then,
	\begin{align*}
	b^\tau_n(S) =  tr \left(Cov\left( \frac{1}{n}\sum_{i=1}^n Z_S(X_{ni})\varphi_\tau\big(Y_{ni} - Z_S(X_{ni})'\theta^\tau_{n,S}\big), \: \hat{\theta}^\tau_{n,S} - \theta^\tau_{n,S} \right) \right) \:\: + \:\: r_{n,1}(S),
	\end{align*}
	where $\varphi_\tau(u) = \tau - 1\{u< 0\}$ and 
	\begin{align*}
	\sup_{S \in M} \left|r_{n,1}(S)\right| = O\left(\frac{1}{\lambda_{n}^{3/2}} \left(\frac{m \log |M| \: \log \log n}{n}\right)^{5/4} \right).
	\end{align*}
\end{theorem}

We postpone a discussion of the rate of the remainder term to the next section. Focusing instead on the leading term of above approximation, we observe the following: If the true CQF is indeed linear in $Z_S(X_n)$, i.e. $Q_{Y_n}(\tau|X_n) = Z_S(X_n)'\theta^\tau_{n,S}$, then the leading term of the optimism $b^\tau_n(S)$ can be re-formulated as
\begin{align}
\frac{1}{n}\sum_{i=1}^n Cov\left(- 1\big\{Y_{ni} < Q_{Y_n}(\tau|X_{ni})\big\}, \:  \widehat{Q}_{Y_n}(\tau|X_{ni}, S) \right).
\end{align}
Thus, in this case the expected optimism is essentially the covariance between the estimates $\widehat{Q}_{Y_n}(\tau|X_{ni}, S)$ and a simple function of the targets $Q_{Y_n}(\tau|X_{ni})$, $i=1, \ldots, n$. This is reminiscent of~\possessivecite{efron2004} results for the ``q class of error measures'': For this specific class of error measures the expected optimism is equal to the covariance between estimate $\hat{\mu}_n$ and target $\mu$, i.e. $\frac{1}{n}\sum_{i=1}^n Cov\big(\mu(Z_i), \: \hat{\mu}_n(Z_i) \big)$. Theorem~\ref{theorem:ExpectOptCovForm} states that for the check loss $\rho_\tau$ a similar relation holds up to the deterministic error $r_{n,1}(S)$.

Re-writing the leading term of the optimism $b^\tau_n(S)$ as the expected value of the gradient of the check loss and the centered regression vector, 
\begin{align}\label{eq:CovarianceFormGradientVersion}
\mathbb{E}_{\mathcal{D}_n}\left[\frac{1}{n}\sum_{i=1}^n\varphi_\tau\big(Y_{ni} - Z_S(X_{ni})'\theta^\tau_{n,S}\big)  Z_S(X_{ni})'(\hat{\theta}^\tau_{n,S}- \theta^\tau_{n,S})\right],
\end{align}
we gain two more insights:

First, the covariance form of the expected optimism can be viewed as a first order linearization of the check loss. In particular, the covariance form is the (expected value) of the directional derivative of the check loss in direction $\hat{\theta}^\tau_{n,S}- \theta^\tau_{n,S}$ and evaluated at the vector of regression coefficients $\theta^\tau_{n,S}$. Since the check loss is convex, this directional derivative is always non-negative, i.e. the leading term of the expected optimism  non-negative. This confirms our statistical intuition that the bias of the in-sample risk as estimate of the predictive risk is negative.

Second, using the naive sample analogue $\frac{1}{n}\sum_{i=1}^n \varphi_\tau\big(Y_{ni} - Z_S(X_{ni})'\hat{\theta}^\tau_{n,S}\big) X_{ni,S}'\hat{\theta}^\tau_{n,S}$ to estimate the expected optimism will inevitably result in a poor estimate because the gradient evaluated at its sample minimizer $\hat{\theta}^\tau_{n,S}$ is close to zero. Thus, even though the approximate covariance form does not dependent on the future (unattainable) data point $(Y^0_n, X^0_n)$, it does not allow us to entirely bypass the computation of the expected value with respect to the unknown distribution $F_n$. A similar observation was first made by~\cite{efron1986} about his covariance penalties. To overcome this difficulty, he proposes a parametric bootstrap approach; below we show a different approach which does not rely on re-sampling.

\subsection{The trace form of the expected optimism}\label{subsec:TraceFormExpectOpt}
As noted in Section~\ref{subsec:PredRiskandExpectOptQR}, the predictive risk under check loss $\rho_\tau$ is an almost affine transformation of the unconditional mean squared prediction error. We might therefore expect that the expected optimism can be approximated by an expression similar to the penalty term in~\possessivecitewos{mallows1973} $C_p$ or~\possessivecite{takeuchi1976} TIC. The following theorem shows that this intuition is correct.

\begin{theorem}[Trace Form of the Expected Optimism]\label{theorem:ExpectOptTraceForm}
	Suppose that Assumptions (A1) -- (A6) from Section~\ref{sec:Assumptions} hold. Then,
	\begin{align*}
	b_n^\tau(S) = \frac{1}{n}tr\left(D_{n,0}^\tau(S)^{-1} D_{n,1}^\tau(S)\right) + r_{n,2}(S),
	\end{align*}
	where
	\begin{align*}
	D_{n,0}^\tau(S) &=\mathbb{E}_{X_{n1}}\Big[f_{Y_n|X_n}\big(Z_S(X_{n1})'\theta^\tau_{n,S}|X_{n1}\big)Z_S(X_{n1})Z_S(X_{n1})'\Big],\\
	D_{n,1}^\tau(S) &=\mathbb{E}_{X_{n1}}\Big[\varphi_\tau^2\big(Y_{n1}- Z_S(X_{n1})'\theta^\tau_{n,S}\big) Z_S(X_{n1})Z_S(X_{n1})'\Big],
	\end{align*}
	with $\varphi_\tau(u) = \tau - 1\{u< 0\}$ and
	\begin{align*}
	\sup_{s \in M} \left|r_{n,2}(S)\right| = O \left(\frac{1}{\lambda_{n}^2}\left(\frac{m \log |M| \: \log\log n}{n}\right)^{5/4} \right).
	\end{align*}
\end{theorem}

We observe the following: First, under Assumptions (A1) -- (A6) the trace from is roughly of order $O\left(\lambda_n^{-1}n^{-1}|S|\right)$ and hence dominates the remainder term $r_{n,2}(S)$. Therefore, the trace form is a meaningful approximation of the expected optimism. We also conclude that the same is true for the covariance form and the remainder term $r_{n,1}(S)$ from Theorem~\ref{theorem:ExpectOptCovForm}.

Second, in the literature on robust estimation the trace form is also known as ``expected self-influence'', i.e. the average influence that an observation has on its own fitted value~\citep[e.g.][p. 317]{hampel2005}. While at hindsight the connection between expected optimism and ``expected self-influence'' appears intuitive, it has not been made in the past, to the best of our knowledge.

Third, the trace form clearly resembles the complexity penalties of AIC-type model selection criteria for misspecified (linear) regression models~\citep[e.g.][]{takeuchi1976, bozdogan2000} and misspecified robust and generalized linear models~\citep[e.g.][]{ronchetti1985, lv2014}. This similarity is expected since complexity penalties of AIC-type model selection criteria aim at estimating the expected optimism of the in-sample risk based on a loss function equal to the negative (pseudo) log-likelihood.

Lastly, by Theorem~\ref{theorem:ExpectOptTraceForm} the expected optimism is a nonlinear function of the conditional density $f_{Y_n|X_n}$, the quantile level $\tau$, the (weighted) covariance of the predictors $Z_S(X_{n})$, and the size $|S|$ of model $S$. This property becomes more salient in the following two special cases:

\begin{corollary}[Location Model]\label{corollary:TraceFormLocationModel}
	Let $Y_{ni} = X_{ni}'\theta_{S_0} + \epsilon_{ni}$, with i.i.d. covariates $X_{ni}$ and i.i.d. errors $\epsilon_{ni} \sim F_\epsilon$ and density $f_\epsilon$. Suppose that the $X_{ni}$ and $\epsilon_{ni}$ are mutually independent for $i=1, \ldots, n$. Let the map $Z$ be the identity map so that the $Z(X_{ni}) = X_{ni}$. Suppose that the conditions of Theorem~\ref{theorem:ExpectOptTraceForm} hold and that the fitted model $S$ contains the true model $S_0$, i.e. $S_0 \subseteq S$. Then,
	\begin{align*}
	\frac{1}{n}tr\left(D_{n,0}^\tau(S)^{-1} D_{n,1}^\tau(S)\right) = \frac{\tau (1- \tau)}{f_{\epsilon}\left(F^{-1}_\epsilon(\tau)\right)}\frac{|S|}{n}.
	\end{align*}
\end{corollary}

\begin{corollary}[Nested Quantile Regression Location Models]\label{corollary:TraceFormNestedLocationQR}
	Suppose that the data generating process is a (potentially nonlinear) location model. Let $S_1$ and $S_2$ be two models such that $S_1 \subseteq S_2$. The trace form of the larger model $S_2$ can be written in terms of the conditional density of $Y_n$ given the predictors $Z_{S_1}(X_{n})$ of the smaller model, i.e.
	\begin{align*}
	\frac{1}{n}tr\left(D_0^\tau(S_2)^{-1} D_1^\tau(S_2)\right) = \frac{\tau (1- \tau)}{n}tr\left(D_0(S_1, S_2)^{-1} D_1(S_2) \right),
	\end{align*}
	where
	\begin{align*}
	D_0(S_1, S_2) &=\mathbb{E}_{X_n}\left[f_{Y_n|Z_{S_1}(X_{n})}\big(Z_{S_2}(X_{n})'\theta^\tau_{S_2}|Z_{S_1}(X_{n})\big)Z_{S_2}(X_{n})Z_{S_2}(X_{n})'\right],\\
	D_1(S_2) &=\mathbb{E}_{X_n}\left[Z_{S_2}(X_{n})Z_{S_2}(X_{n})'\right].
	\end{align*}
\end{corollary}
Both corollaries are an immediate consequences of Theorem~\ref{theorem:ExpectOptTraceForm} and~\possessivecite{angrist2006} characterization of the misspecified quantile regression problem as a weighted least squares problem. We omit their proofs. We will return to these two corollaries in Section~\ref{sec:NumericalEvidence} and use them as benchmark in our numerical experiments.

\section{Consistent estimators for expected optimism and predictive risk}\label{sec:EstimationPredRisk}
\subsection{A plug-in estimator for the expected optimism}\label{subsec:EstimationExpectOpt}
The trace form of Theorem~\ref{theorem:ExpectOptTraceForm} lends itself to a simple plug-in estimator for the expected optimism since the two matrices $D_{n,0}^\tau(S)$ and $D_{n,1}^\tau(S)$ are well-studied in the context of the (asymptotic) covariance matrix of the quantile regression vector~\citep[e.g.][]{koenker2005}. In the case of incorectly specified quantile regression models, the following estimates for $D_{n,0}^\tau(S)$ and $D_{n,1}^\tau(S)$ have been proposed
\begin{align}
\widehat{D}^{\tau}_{0,h}(S)  &= \frac{1}{2nh} \sum_{i=1}^{n}1\big\{\big|Y_{ni} - \widehat{Q}_{Y_n}(\tau|X_{ni}, S)\big|\leq h\big\}Z_S(X_{ni})Z_S(X_{ni})',\\
\widehat{D}^\tau_{n,1}(s) &= \frac{1}{n} \sum_{i=1}^{n}\varphi_\tau\big(Y_{ni} - \widehat{Q}_{Y_n}(\tau|X_{ni},S)\big) Z_S(X_{ni})Z_S(X_{ni})',
\end{align}
where $h$ is a bandwidth parameter and $\varphi_\tau(u) = \tau - 1\{u < 0\}$~\citep[e.g.][]{angrist2006, belloni2016}. We therefore propose the following plug-in estimate for the expected optimism $b^\tau_n(S)$,
\begin{align}\label{eq:EstimateExpectOpt}
\hat{b}_{n, h}^\tau(S) = \frac{1}{n}tr\left(\widehat{D}^{\tau^{-1}}_{0,h}(S) \widehat{D}^\tau_{n,1}(S)\right).
\end{align}

Since our regularity conditions are slightly more general than those in~\cite{belloni2016}, the following consistency theorem does not follow from their Lemma 30. In particular, our Assumption (A5) on the growth rate of the number of predictors is less stringent than theirs. We shall see that this relaxation is important in the context of predictive risk estimation in Section~\ref{subsec:EstimationPredRisk}.

\begin{proposition}[Uniform Consistency of the Estimated Trace Form]\label{proposition:ConsistencyEstTraceForm}
	Suppose that Assumptions (A1) -- (A6) from Section~\ref{sec:Assumptions} hold, let $h > 0$ be the bandwidth parameter, and $r_{n} = \frac{1}{\lambda_n}\left(\frac{m \log |M|\: \log \log n}{n}\right)^{1/2}$. Then,
	\begin{align*}
	\sup_{S \in M} \left|n \cdot \hat{b}_{n, h}^\tau(S) - tr\Big(D_{n,0}^\tau(S)^{-1} D_{n,1}^\tau(S)\Big) \right| &= O_p \left( \frac{m \: h^\alpha}{\lambda_n^2} + \frac{m\:r_{n}}{h \lambda_n} + \frac{m\: r_{n}^\alpha}{\lambda_n^2} \right).
	\end{align*}
\end{proposition}
The first and second terms on the right hand side capture the variance and bias of the estimator with bandwidth $h$. They are standard in nonparametric smoothing. The third term controls the bias induced by $\big\{\big(Y_{ni} - \widehat{Q}_{Y_n}(\tau|X_{ni}, S)\big)\big\}_{i=1}^n$ at model $S$ which serve as proxies for $\big\{\big(Y_n - Z_S(X_{ni})'\theta_{n,S}^\tau\big)\big\}_{i=1}^n$.

Specializing to the common case of a continuous conditional density $f_{Y_n|X_n}$, i.e. $\alpha = 1$, we observe the following: The optimal, mean-variance-balancing bandwidth is $h^* = (c_1/c_0)^{1/2}(\lambda_n r_{n})^{1/2}$ with constants $c_0, c_1 > 0$ given in eq.~\eqref{eq:BoundAA1} and~\eqref{eq:BoundAA2}, respectively. In principle, these constants can be estimated from the data; however, in practice, we find that the specific choice of the bandwidth has no significant effect. With bandwidth $h^*$ the estimate $\hat{b}_{n, h}^\tau(S)$ is consistent at rate $O_p\big(m \: r_{n}^{1/2} \lambda_n^{-3/2} + m \: r_{n} \lambda_n^{-2} \big) = O_p\big(m \: r_{n}^{1/2} \lambda_n^{-3/2}\big)$. That is, $\hat{b}_{n, h}^\tau(S)$ is consistent at a rate that is the same as if the true errors $\big\{\big(Y_{ni} - \widehat{Q}_{Y_n}(\tau|X_{ni}, S)\big)\big\}_{i=1}^n$ at model $S$ were known.

Combining Theorem~\ref{theorem:ExpectOptTraceForm} and Proposition~\ref{proposition:ConsistencyEstTraceForm} we obtain the following consistency result.
\begin{theorem}[Uniform Consistency of the Estimated Expected Optimism]\label{theorem:ConsistencyEstExpectOpt} Let $r_{n} = \frac{1}{\lambda_n}\left(\frac{m \log |M|\: \log \log n}{n}\right)^{1/2}$. Under the conditions of Proposition~\ref{proposition:ConsistencyEstTraceForm},
	\begin{align*}
	\sup_{S \in M}\left| \frac{\hat{b}_{n, h}^\tau(S)}{b_n^\tau(S)} -1 \right|  &= O_p\left( n \lambda_n^{3/2} r_{n}^{5/2} +  \frac{m\: h^\alpha}{\lambda_n^2} + \frac{m\: r_{n}}{h\lambda_n} + \frac{m \: r_{n}^\alpha}{\lambda_n^2}\right).
	\end{align*}
\end{theorem}
Since $\hat{b}_{n, h}^\tau(S)$ is the plug-in estimator for the trace form approximation, it is a biased estimate of the actual expected optimism $b_n^\tau(S)$. This deterministic bias is captured in the first term; the remaining three terms are already familiar from Proposition~\ref{proposition:ConsistencyEstTraceForm}. Specializing once again to the case of a continuous conditional density, i.e. $\alpha = 1$, we have under the optimal bandwidth $h^*$ a rate of $O_p\big(n \lambda_n^{3/2} r_{n}^{5/2} + m \: r_{n}^{1/2} \lambda_n^{-3/2}\big) = O_p\big(n \lambda_n^{3/2} r_{n}^{5/2}\big)$. Thus, the deterministic error of using the trace form $tr\big(D_{n,0}^\tau(S)^{-1} D_{n,1}^\tau(S)\big)$ to approximate the expected optimism $b_n^\tau(S)$ dominates the stochastic estimation error. In other words, as point estimate $\hat{b}_n^\tau(S)$ is as good in estimating the expected optimism $b_n^\tau(S)$ as the unattainable trace form $tr\big(D_{n,0}^\tau(S)^{-1} D_{n,1}^\tau(S)\big)$.

\subsection{A de-biased estimator of the predictive risk}\label{subsec:EstimationPredRisk}
As outlined in Section~\ref{subsec:PredRiskandExpectOptGENERAL}, given the consistent estimate of the expected optimism~\eqref{eq:EstimateExpectOpt} we can construct the following de-biased estimate of the predictive risk,
\begin{align}\label{eq:EstimatePredRisk}
\widehat{PR}^{\tau}_{n,h}(S) = \frac{1}{n} \sum_{i=1}^n \Big(\rho_\tau\big(Y_{ni} - \widehat{Q}_{Y_{n}}(\tau|X_{ni}, S)\big) - \rho_\tau(Y_{ni})\Big) + \hat{b}_{n, h}^\tau(S).
\end{align}
We call this estimate ``de-biased'' because the in-sample risk $\frac{1}{n} \sum_{i=1}^n \big(\rho_\tau\big(Y_{ni} - \widehat{Q}_{Y_{n}}(\tau|X_{ni}, S) - \rho_\tau(Y_{ni})\big)$ is itself already a consistent estimate for $\mathrm{PR}^{\tau}_{n,h}(S)$ in the sense that for any $S \in M$ with fixed model size $|S|$,
\begin{align}\label{eq:PredRiskNaive}
\left|\mathrm{PR}^{\tau}_{n}(S) - \frac{1}{n} \sum_{i=1}^n \Big(\rho_\tau\big(Y_{ni} - \widehat{Q}_{Y_{n}}(\tau|X_{ni}, S)\big) - \rho_\tau(Y_{ni})\Big)\right| = O_p\left(n^{-1/2}\right).
\end{align}

We strengthen this fact in several ways: First, we show that under appropriate conditions our proposed estimator $\widehat{PR}^{\tau}_{n,h}(S)$ is consistent uniformly over all $S \in M$ and for models whose size $|S|$ grows with the sample size $n$. Second, we will see that for large models with size $|S| \gtrsim n^{1/2}$ the in-sample risk is no longer $n^{1/2}$-consistent for the predictive risk and that under certain conditions de-biasing the in-sample risk with $\hat{b}_{n, h}^\tau(S)$ restores the $n^{1/2}$-consistency. We deduce these claims from the following general result.

\begin{theorem}[Uniform Consistency of the De-biased Predictive Risk Estimate]\label{theorem:PredRiskConsistency}
	Suppose that Assumptions (A1) -- (A6) from Section~\ref{sec:Assumptions} hold. In addition, assume that $f_{Y_{n}|X_{n}}$ is uniformly bounded away from 0 for all $n$ and that $\limsup_{n \rightarrow \infty}\mathbb{E}_{X_n}\left[Q_{Y_n}^2(\tau|X_n)\right] < \infty$.  Let $h > 0$ be a bandwidth and $r_n = \frac{1}{\lambda_n}\left(\frac{m \log|M| \: \log\log n}{n}\right)^{1/2}$. Then,
	\begin{align*}
	& \sup_{S \in M}\left| \widehat{PR}^{\tau}_{n,h}(S) - \mathrm{PR}^\tau_n(S) \right| = O_p \left(\left( \frac{\log |M|}{n}\right)^{1/2} + \frac{r_n}{n^{1/2}} + \lambda_n^{3/2} r_{n}^{5/2} + \frac{m \: h^\alpha}{\lambda_n^2 n} + \frac{m\:r_{n}}{h \lambda_n n} + \frac{m\: r_{n}^\alpha}{\lambda_n^2 n}\right),
	\end{align*}
\end{theorem}

The last four terms on the right hand side are familiar from the uniform consistency result of the trace form estimate for the expected optimism (i.e. Theorem~\ref{theorem:ConsistencyEstExpectOpt}), while the first two terms are related to the in-sample risk. Clearly, if $m= o( \lambda_n^{-2} n \log |M| \: \log \log n)$ and bandwidth $h$ satisfies
\begin{align}
\frac{1}{\lambda_n} \frac{m}{n} \left(\frac{m \log|M| \: \log\log n}{n}\right)^{1/2} \lesssim h \lesssim \frac{1}{\lambda_n^{2/\alpha}} \left(\frac{n}{m}\right)^{1/\alpha},
\end{align}
then $\widehat{PR}^{\tau}_{n,h}(S)$ is consistent for $\mathrm{PR}^\tau_n(S)$ uniformly for all $S \in M$. However, we can learn more by considering special cases. To simplify this discussion, we consider the case in which the conditional density $f_{Y_n|X_n}$ is continuous and the bandwidth is chosen to balance the nonparametric estimation bias and variance (see discussion in Section~\ref{subsec:EstimationExpectOpt}). Then, Theorem~\ref{theorem:PredRiskConsistency} implies the following.

\begin{corollary}\label{corollary:PredRiskConsistencyContinuous}
	Suppose that the conditions of Theorem~\ref{theorem:PredRiskConsistency} hold, that the conditional density $f_{Y_n|X_n}$ is continuous, and that $ \lambda_n^2 m = o\left(n \log |M| \: \log \log n\right)$. If $ n^{1/4} h \sim \left(m \log|M| \: \log\log n\right)^{1/4}$, then
	\begin{align*}
	& \sup_{S \in M}\left| \widehat{PR}^{\tau}_{n,h}(S) - \mathrm{PR}^\tau_n(S) \right| = O_p \left(\left( \frac{\log |M|}{n}\right)^{1/2} + \frac{1}{\lambda_n^2} \left(\frac{m \log|M| \: \log\log n}{n}\right)^{5/4} \right).
	\end{align*}
\end{corollary}

These rates have an intuitive explanation: The first term $O\big(n^{-1/2}(\log |M|)^{1/2}\big)$ is related to the stochastic variability of the in-sample risk. And the second term $O(\lambda_n^{-2} n^{-5/4} (m \log|M| \log \log n)^{5/4})$ is known from Theorem~\ref{theorem:ExpectOptTraceForm} to be the deterministic error of using the trace form $tr(D_{n,0}^\tau(S)^{-1} D_{n,1}^\tau(S))$ to approximate the expected optimism $b_n^\tau(S)$. Thus,  unlike one  might have suspected, it is not the nonparametric estimate of the expected optimism but the deterministic approximation of the expected optimism and the stochastic variability of the in-sample risk which limit the accuracy of our predictive risk estimate. It is easy to verify that under the stated assumptions $\widehat{PR}^{\tau}_{n,h}(S)$ is consistent for $\mathrm{PR}^\tau_n(S)$ uniformly over all $S \in M$. 

It is instructive to consider the implication of Corollary~\ref{corollary:PredRiskConsistencyContinuous} under different growth regimes of the number of predictor variables. To this end, recall that the estimated trace form, $\hat{b}_{n, h}^\tau(S)$, is of order $O(\lambda_n^{-1}n^{-1}|S|)$.
Hence, if $n^{1/2} \lesssim |S| \lesssim n$ the estimated trace form, $\hat{b}_{n, h}^\tau(S)$, dominates (rate-wise) the stochastic error and also the deterministic error (provided that we sharpen condition on $m$ and $n$ to $m = o(n \: \lambda_n^4 (\log |M| \: \log \log n)^{-5})$). Thus, in this regime the in-sample risk alone is not $n^{1/2}$-consistent for the predictive risk; de-biasing the in-sample risk is necessary to retain $n^{1/2}$-consistency.

However, if $|S| \lesssim n^{1/2}$ the stochastic error of the in-sample risk dominates (rate-wise) the estimate of the trace form. Thus, from the perspective of first order asymptotics the correction provided by the $\hat{b}^\tau_{n,h}(S)$ is not necessary in this regime. However, in Section~\ref{sec:NumericalEvidence} we report numerical evidence showing that even in this regime the de-biasing effect of $\hat{b}^\tau_{n,h}(S)$ is practically relevant.

As an aside, this discussion provides another explanation for the well-known fact that Akaike-type model selection criteria are not model selection consistent: Akaike-type penalties (based on estimates of the expected optimism) are too small to effectively discriminate between models of size $|S| \lesssim n^{1/2}$ since the stochastic variability of the in-sample risk is relatively large. For correctly specified (linear least squares regression) models with a fixed number of parameters this has already been recognized~\citep[e.g.][]{shao1997, yang2005}.

It is natural to consider using the de-biased predictive risk estimator for model selection purposes. Indeed, the uniform consistency result from Theorem~\ref{theorem:PredRiskConsistency} implies that the model(s) minimizing the de-biased predictive risk estimate $\widehat{PR}^{\tau}_{n,h}$ are consistent for the model(s) minimizing the predictive risk $\mathrm{PR}^\tau_n$. For a precise statement and proof of this claim, we refer to the Supplementary Materials~\ref{Supplement:ModelSelection}. In contrast, neither Akaike-type Information Criteria~\citep[e.g.][]{burman1995, koenker2005}, nor Schwarz/ Bayesian Information Criteria~\citep[e.g.][]{machado1993, koenker1994, lee2014}, nor the Asymptotically defined Information Criterion~\citep{portnoy1997} are known to satisfy such a model selection consistency property if the candidate models are (possibly) misspecified.

\section{Empirical evidence}\label{sec:NumericalEvidence}
\subsection{Set-up of the simulation study}
We conduct Monte Carlo experiments to evaluate empirically the trace form approximation of the expected optimism and to corroborate the theoretical results from Sections~\ref{sec:ApproxExpectOpt} and~\ref{sec:EstimationPredRisk}. We also compare the empirical performance of the trace form approximation to the commonly used cross-validated estimate of the expected optimism. Our Monte Carlo study uses four designs as the data generating processes (DGP), but only the results from DGP1 are given in the paper. The results from the other DGPs are qualitatively similar and details are given in the Supplementary Materials in~\ref{Supplement:NumericalEvidence}.

\begin{itemize}
	\item Independent Gaussian design (DGP1): $y_i = x_{i1} + x_{i2} + x_{i3} + x_{i4} + \epsilon_i$, with $x_{i} \sim_{iid} N(0, I_p)$ independent of the errors $\epsilon_i \sim_{i.i.d.} N(0, 4)$.
	
	We use this process to illustrate the elementary properties of the predictive risk and the expected optimism from Corollaries~\ref{corollary:TraceFormLocationModel} and~\ref{corollary:TraceFormNestedLocationQR}. The joint Gaussianity of predictors and errors allows us to compute the exact value of the trace form with which we can assess the accuracy of our estimates. The variance of the error distribution is chosen such that signal-to-noise-ratio equals one.
	
	\item Correlated Gaussian design (DGP2): $y_i = x_{i1} + x_{i2} + x_{i3} + x_{i4} + \epsilon_i$, with  $\epsilon_i \sim_{i.i.d.} N(0, 12.384)$ independent of $x_{i} \sim_{iid} N(0, \Sigma)$ and $\Sigma_{ij}= 0.8^{|i-j|}$ for all $i, j = 1, \ldots, p$.

	The variance of the error distribution is chosen such that the signal-to-noise ratio equals one.

	\item Heteroscedastic noise (DGP3): $y_i = x_{i1} + x_{i2} + x_{i3} + (1 + 1.5 x_{i4})\epsilon_i$, where $x_{ij} \sim_{i.i.d.} U\left([0,2]\right)$ for $j=1, \ldots, 4$ independent of the errors $\epsilon_{i} \sim_{iid} N(0, 1)$.

	In this DGP the covariate $x_4$ is active for the conditional quantile functions except at the median.

	\item Single interaction term with heavy-tailed noise (DGP4): $y_i = x_{i1} + x_{i2} + x_{i3} + 4x_{i3}x_{i4} + \epsilon_i$, where $\epsilon_i$ follow the t-distribution with 2 degrees of freedom independent of the predictors $x_{i} \sim_{iid} N(0, I_p)$.
	
	In this DGP all quantiles are non-linear functions of the covariates.
\end{itemize}

We set the dimension of the space of covariates $\mathcal{X}$ equal to 50, and let $Z$ be the identity map, so that the predictors are simply the covariates $X_1, \ldots, X_{50}$. We consider a collection of 176 candidate models with model seizes ranging between 0 to 50. This implies that we the size of the largest model under consideration is $m=50$. We explain the choice of those candidate models in Section~\ref{subsec:TraceFormIS}. Throughout the numerical experiments we keep the sample size fixed at $n=500$. All reported estimates are averages over 10,000 independent realizations of the corresponding DGPs. To estimate the matrix $D_0(S)$ at quantile $\tau$ we use~\possessivecite{powell1986} nonparametric estimator with uniform kernel function and bandwidth
\begin{align*}
c_{n,S} = \kappa_{n,S} \left(\Phi^{-1}(\tau + h_n) - \Phi^{-1}(\tau - h_n) \right),
\end{align*}
where $\Phi$ denotes c.d.f. of the standard normal distribution, $\kappa_{n,S}$ is the minimum of the standard error and the inter-quartile-range of the estimated quantile regression residuals of model $S$, and
\begin{align*}
h_n = \frac{1}{n^{1/5}}\left(\frac{4.5\phi\big(\Phi^{-1}(\tau)\big)^4}{(2\Phi^{-1}(\tau)^2 + 1)^2}\right)^{1/5},
\end{align*}
where $\phi$ denotes the p.d.f. of the standard normal distribution. Thus, $c_{n,S}$ satisfies the conditions of Theorems~\ref{theorem:ConsistencyEstExpectOpt} and~\ref{theorem:PredRiskConsistency} which guarantee (uniform) consistency of the estimates of the expected optimism and the predictive risk; see~\cite{koenker2005} for a detailed discussion of this choice of bandwidth.

Recall Definitions~\ref{def:PredRisk} and~\ref{def:ExpectedOptimism} that the predictive risk and the expected optimism require the evaluation of a double expectation. Since the quantile regression vector is only implicitly defined, this double expectation cannot be evaluated analytically. Instead, we use Monte Carlo estimates based on 50,000 samples to obtain values for the predictive risk and the expected optimism.

\subsection{Estimation of the expected optimism}\label{subsec:TraceFormIS}
\begin{figure}[!htbp]	\centering
	\minipage{0.85\textwidth}
	\includegraphics[width=\linewidth]{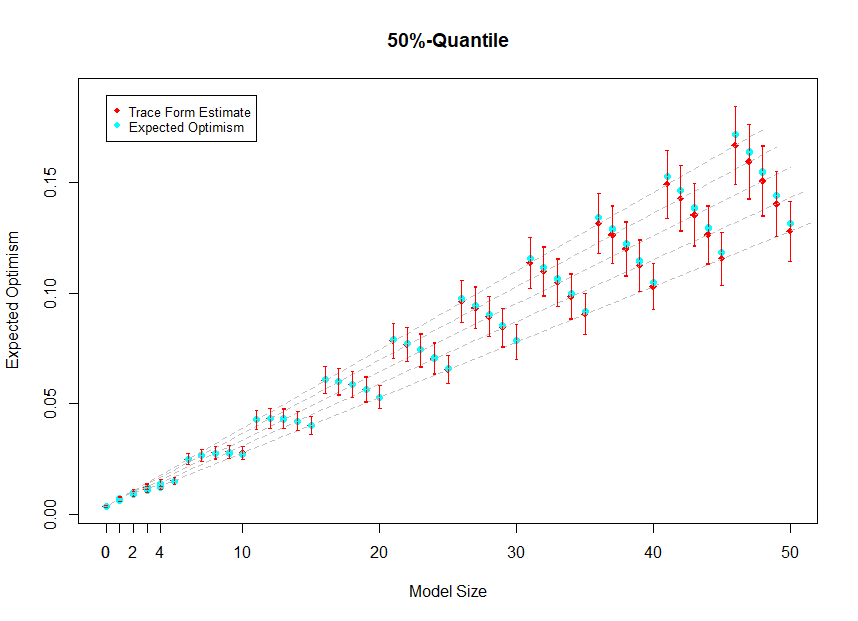}
	\endminipage\hfill\\
	\minipage{0.85\textwidth}
	\includegraphics[width=\linewidth]{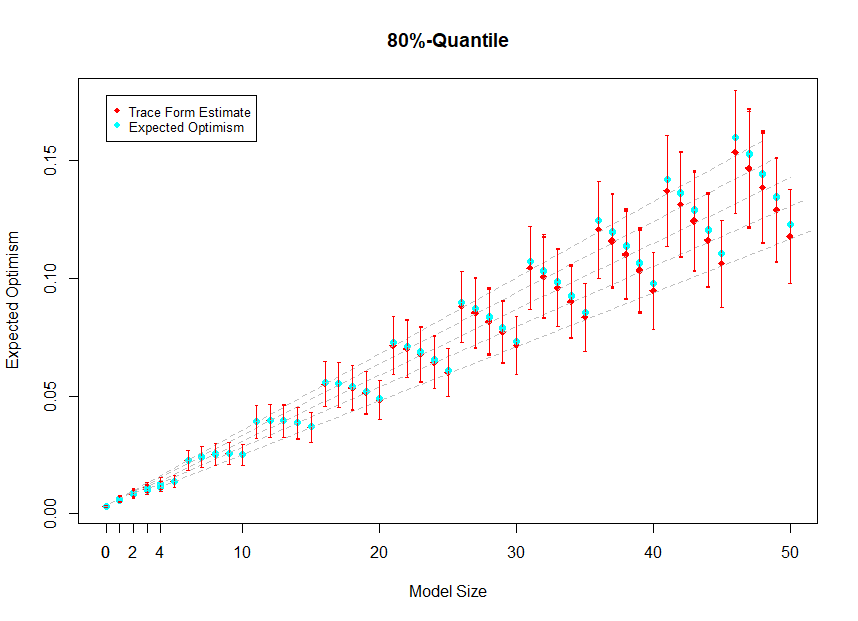}
	\endminipage
	\caption{Red: estimates of the trace form and standard errors. Blue: expected optimism. Dashed gray lines: exact evaluation of the trace form. Top: DGP1 with $\tau = 0.5$. Bottom: DGP1 with $\tau =0.8$. }\label{fig:DGP1-IS}
\end{figure}

In Theorem~\ref{theorem:ConsistencyEstExpectOpt} we establish uniform consistency of the estimated trace form for the expected optimism.
In Figure~\ref{fig:DGP1-IS} under DPG1, we plot the bias of 176 models (subsets of the 50 predictors) against their model sizes. We only consider 176 models because it is computationally expensive to evaluate the predictive risk and the expected optimism on all possible subsets of the 50 predictors. However, the special structure of the DGP together with Corollary~\ref{corollary:TraceFormNestedLocationQR} guarantee that this collection constitutes a representative subset of all possible models: The true DGP contains only four relevant predictors 1, 2, 3, and 4; those predictors are independent and identically distributed and contribute equally to the model (i.e. have the same regression coefficients). We can therefore stratify the collection of all possible subsets of the 50 predictors according to how many relevant predictors are included in a specific subset. This results in five collections of nested models indexed by 0 (relevant predictors), 1 (relevant predictor), \ldots, 4 (relevant predictors). By Corollary~\ref{corollary:TraceFormNestedLocationQR} the expected optimism of all nested models with $j$ relevant predictors lie (approximately) on a ray emanating from the in-sample bias of the smallest model with $j$ relevant predictors. Moreover, the slope of the ray is given by $\frac{\tau(1-\tau)}{500 \phi_j}$, where $\phi_j$ denotes the value of the density of a centered normal random variable with variance $j^2 + 1$ evaluated at 0. The 176 models comprise the model that contains only the intercept and 35 models of each of the five stratified collections. 

In Figure~\ref{fig:DGP1-IS} the top gray line corresponds to the theoretical values of the trace form of models that have four relevant predictors and additional, irrelevant, predictors. The second line from the top corresponds to the theoretical values of the trace form of models that contain three relevant predictors and additional, irrelevant, predictors, and so forth. The last line (fifth from above) corresponds to models that do not contain any relevant predictors.

We observe that the estimates of the trace form (in red) lie on (or are very close) to theoretical values of the trace form uniformly for all 176 models. This confirms the fast uniform convergence rates obtained in Theorem~\ref{theorem:ConsistencyEstExpectOpt}. Note that the plot shows only 50 red dots and not as one might expect 176 dots. This is due to the fact that for DGP1 the value of the estimated trace form does not depend on the specific subset of predictors (i.e. $S$) but only on the size of the model (i.e. $|S|$), e.g. the two models with predictors $\{1,2, 5\}$ and $\{3,4, 10\}$ have the same trace form which is fully determined by the fact that they contain two relevant and one irrelevant predictors. The expected optimism (in blue) does not follow the dashed gray lines of the theoretical values of the trace form as closely as the estimates do. This reflects the fact that the trace form is only an approximation to the expected optimism (see Theorem~\ref{theorem:ExpectOptTraceForm}). The difference between the values of the trace form and the expected optimism appears to be negligible for models of size up to $20 \approx \sqrt{n}$ (recall that $n= 500$).

The vertical red lines indicate the standard deviations of the estimated trace forms. The standard deviation increases with the model size and, holding the number of nuisance predictor variables fixed, decreases with the number of relevant predictor variables that are included in the model. The latter effect is rather weak and can be best observed in the plot for the 80\% quantile.

\subsection{Comparison with cross-validated expected optimism}
Cross-validation is a commonly-used method for estimating the predictive risk and the expected optimism.
In this subsection we compare the trace form estimate with a 10-fold cross-validation estimate of the expected optimism.

\begin{figure}[!htbp]
	\minipage{0.5\textwidth}
	\includegraphics[width=\linewidth]{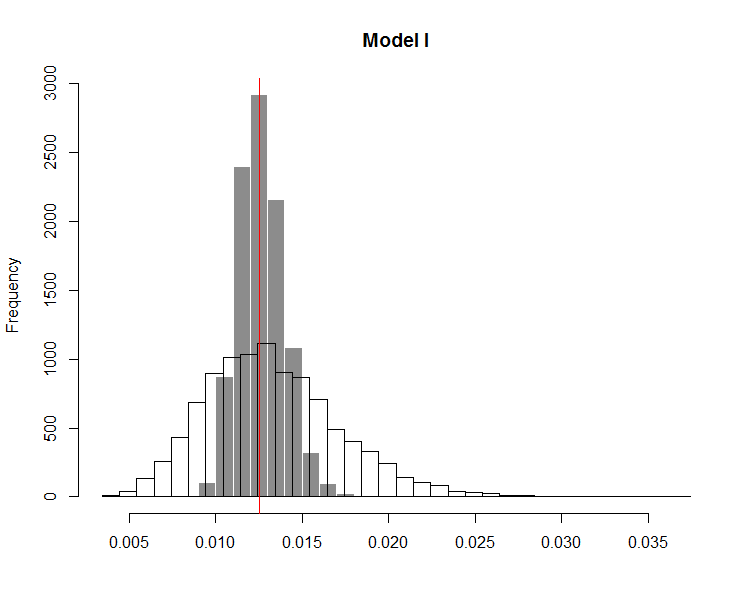}
	\endminipage\hfill
	\minipage{0.5\textwidth}
	\includegraphics[width=\linewidth]{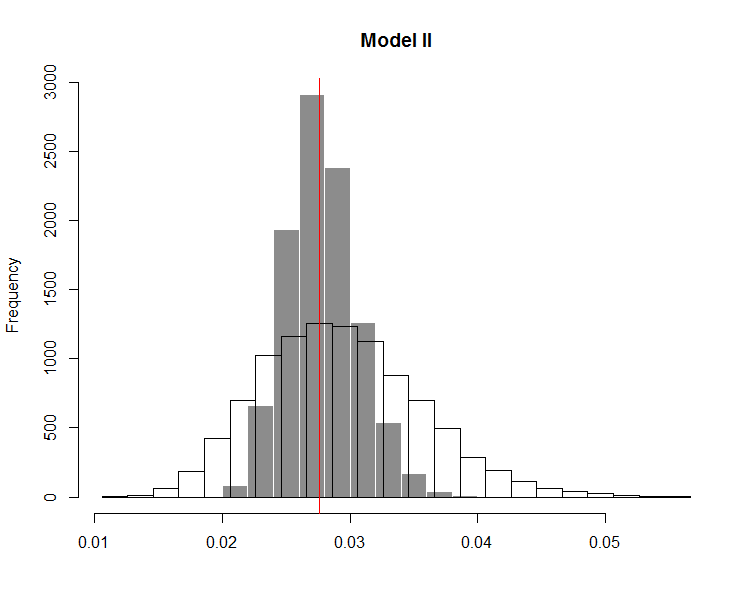}
	\endminipage\hfill
	\minipage{0.5\textwidth}
	\includegraphics[width=\linewidth]{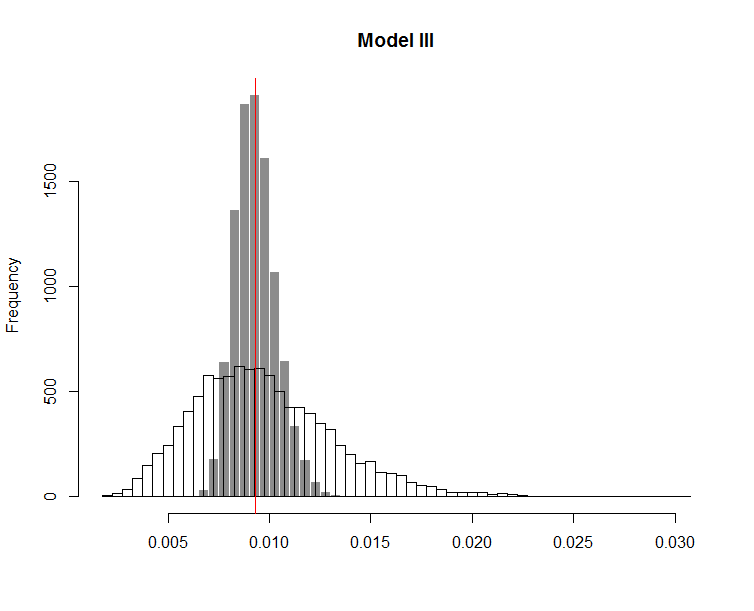}
	\endminipage\hfill
	\minipage{0.5\textwidth}
	\includegraphics[width=\linewidth]{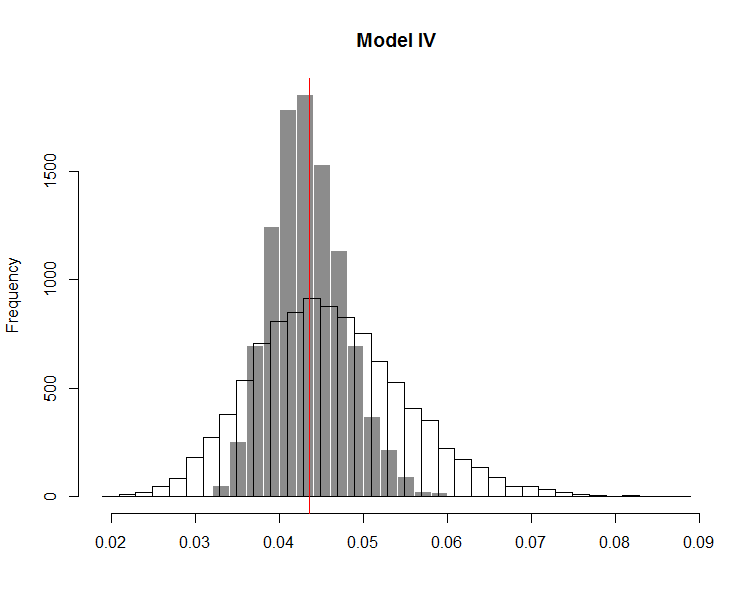}
	\endminipage
	\caption{Histograms of the 10-fold CV estimate of the expected optimism and the trace form estimate for DGP1 and $\tau = 0.5$. Red line: expected optimism. White histogram: 10-fold CV. Gray histogram: trace form estimate. Model I: correct model (with predictors 1 to 4), Model II: an over-fitted model (with predictors 1 to 10), Model III: an under-fitted model (with predictors 1 to 2) and Model IV that comprises the relevant predictors 1 and 2 and the irrelevant predictors 5 to 15.}\label{fig:DGP1-IS-Bias}
\end{figure}

Figure~\ref{fig:DGP1-IS-Bias} shows the results of 10-fold cross-validation and the trace form for DGP1 at the median. We consider four representative models: Model I is the correct model (with predictors 1 to 4), Model II is an over-fitted model (with predictors 1 to 10), Model III is an under-fitted model (with predictors 1 and 2) and Model IV is the model that comprises the relevant predictors 1 and 2 and the irrelevant predictors 5 to 15. The vertical red line indicates the expected optimism. The white histograms show the empirical distribution of 10,000 cross-validation estimates of the expected optimism and the dark gray histograms show the empirical distribution of 10,000 trace form estimates of the expected optimism.

Both histograms are centered around the expected optimism; however, the estimate of the trace form concentrates significantly more around the target. As mentioned in Section~\ref{subsec:PredRiskandExpectOptGENERAL} the reason for this is that the cross-validation estimate is based on a smaller sample size both for model estimation and for risk estimation.

\section{Conclusion}\label{sec:Conclusion}
In the present paper, we have derived two asymptotic approximations of the expected optimism, or the bias of the in-sample risk when used as an estimate of the predictive risk, and have proposed consistent estimates of the expected optimism and the predictive risk of potentially misspecified quantile regression models. The asymptotic approximations based on two explicit forms help us understand how the expected optimism depends on several factors, including the quantile level, the model misspecification bias, the model size, and sampling variability. In some simpler cases, the expected optimism is asymptotically linear in the model size, but for under-fitted or misspecified models in general, the relationship is far more complicated. The results show that commonly used AIC-type model selection criteria for quantile regression are not really good proxies of the predictive risk. We propose a bias-corrected estimate of the predictive risk and establish its uniform consistency under weak assumptions. Our theoretical results indicate that de-biasing the in-sample risk with an estimate of the expected optimism is necessary when considering models whose dimension grow with at least $n^{1/2}$. Empirical evidence suggests that even in the case of models with fixed dimension estimates of the predictive risk can be significantly improved via de-biasing the in-sample risk.

The asymptotic approximations derived in the present paper are uniform in a class of candidate models, but those models are not data-dependent. An interesting question that relates more to  model selection criteria is how well the bias, and thus the predictive risk estimation, hold up for data-dependent models. Clearly, additional research is needed to address this question. 

\section*{Acknowledgements}
The research is partly supported by NSF Award DMS-1607840 (USA) and Natural Science Foundation of China Grant 11690012. We thank the Editors and two Referees for their comments and suggestions that have led to significant improvements in presentation and scope of the paper. 

\bibliography{research2}

\appendix
\section{Proofs}\label{appendix:Proofs}
We denote the check loss of the $\tau$th quantile by $\rho_\tau(u) = u(\tau - 1\{u < 0\})$ and the corresponding score function by $\varphi_\tau(u) = \tau - 1\{u < 0\}$. We define $Z_{ni,S} = Z_S(X_{ni})$, $Z^0_{ni,S} = Z_S(X^0_{ni})$, $\hat{\delta}^\tau_{n,S} = \hat{\theta}^\tau_{n,S} - \theta^\tau_{n,S}$, $e^\tau_{ni,S} = Y_{ni} - Z_{ni,S}'\theta^\tau_{n,S}$, $\hat{e}^\tau_{ni,S} = Y_{ni} -Z_{ni,S}'\hat{\theta}^\tau_{n,S}$, and $\hat{e}^{0\tau}_{n,S} = Y^0_{n} -Z_{n,S}^{0'}\hat{\theta}^\tau_{n,S}$. We use $C, c, c_0, c_1, \ldots$ to denote absolute constant that may change from line to line. Let $r_n = \frac{c_3}{\lambda_n}\left(\frac{m \log |M| \:\log \log n}{n} \right)^{1/2}$, where $c_3 > 0$ is the absolute constant from Lemma~\ref{lemma:StrongBahadurRep}. Throughout we assume that $|M| \geq 2$ and $\log \log n > 1$, i.e. $n > 15$.

\subsection{Proof of Theorem~\ref{theorem:ExpectOptCovForm}}
\begin{proof}[\textbf{Proof of Theorem~\ref{theorem:ExpectOptCovForm}}]
	
	\textbf{Step 1:} By Knight's identity,
	\begin{align*}
	\rho_\tau(u-v) -\rho_\tau(u) = - v\varphi_\tau(u) + \int_0^v \left(1\{u \leq s\} - 1\{u\leq 0\}\right)ds,
	\end{align*}
	for arbitrary $S \in M$, we can write the optimism as
	\begin{align}\label{eq:KnightsExpansion}
	&\mathbb{E}_{\mathcal{D}_n, (Y^0_n,X^0_n)}\left[\frac{1}{n}\sum_{i=1}^n \rho_\tau(Y_{ni} - Z_{ni,S}'\hat{\theta}_{n, s}^\tau) - \rho(Y_{ni}) - \rho_\tau(Y^0_{n} - Z_{n,S}^{0'}\hat{\theta}_{n, S}^\tau) + \rho_\tau(Y^0_{n})\right] \nonumber\\
	&= \mathbb{E}_{\mathcal{D}_n}\left[-\frac{1}{n}\sum_{i=1}^n Z_{ni,S}'\hat{\delta}^\tau_{n,S} \varphi_\tau(e^\tau_{ni,S}) + \frac{1}{n}\sum_{i=1}^n \int_0^{Z_{ni,S}'\hat{\delta}^\tau_{n,S}} \left(1\{e^\tau_{ni,S} \leq t\} - 1\{e^\tau_{ni,S}\leq 0\}\right)dt \right] \nonumber\\
	&\quad{} - \mathbb{E}_{\mathcal{D}_n, (Y^0_n, X^0_n)}\left[- Z_{n,S}^{0'}\hat{\delta}^\tau_{n,S} \varphi_\tau(e^{0\tau}_{n,S}) + \int_0^{Z_{n,S}^{0'}\hat{\delta}^\tau_{n,S}} \left(1\{\hat{e}^{0\tau}_{n,S}\leq t\} - 1\{\hat{e}^{0\tau}_{n,S}\leq 0\}\right)dt  \right] \nonumber\\
	& = \left(-\mathbb{E}_{\mathcal{D}_n}\left[\frac{1}{n}\sum_{i=1}^n Z_{ni,S}'\hat{\delta}^\tau_{n,S} \varphi_\tau(e^\tau_{ni,S})\right]  + \mathbb{E}_{\mathcal{D}_n, (Y^0_n, X^0_n)}\left[Z_{n,S}^{0'}\hat{\delta}^\tau_{n,S} \varphi_\tau(\hat{e}^{0\tau}_{n,S})\right]\right) \nonumber\\
	&\quad{} + \left(\mathbb{E}_{\mathcal{D}_n}\left[\frac{1}{n}\sum_{i=1}^n \int_0^{Z_{ni,S}'\hat{\delta}^\tau_{n,S}} \left(1\{e^\tau_{ni,S} \leq t\} - 1\{e^\tau_{ni,S}\leq 0\}\right)dt \right] \right. \nonumber\\
	&\quad{} \left. - \mathbb{E}_{\mathcal{D}_n, (Y^0, X^0)}\left[\int_0^{Z_{n,S}^{0'}\hat{\delta}^\tau_{n,S}} \left(1\{\hat{e}^{0\tau}_{n,S} \leq t\} - 1\{\hat{e}^{0\tau}_{n,S}\leq 0\}\right)dt  \right]\right) \nonumber\\
	& = A_n(S) + B_n(S).
	\end{align}
	
	\textbf{Step 2: Uniform upper bound on $B_n(S)$.} Let $\varepsilon_1, \ldots, \varepsilon_n$ be independent Rademacher random variables. Then,
	\begin{align}
	\sup_{S \in M} B_n(S) & \leq  \sup_{S \in M} \mathbb{E}_{\mathcal{D}_n, (Y^0_n, X^0_n)}\left[\sup_{\|\delta_S\|_2 \leq r_n} \frac{1}{n}\sum_{i=1}^n \int_0^{Z_{ni,S}'\delta_S} \Big(1\{e^\tau_{ni,S} \leq t\} - 1\{e^\tau_{ni,S}\leq 0\}\Big)dt \right. \nonumber\\
	& \left. \phantom{\sup_{S \in M} \mathbb{E}_{\mathcal{D}_n, (Y^0_n, X^0_n)}\left[\sup_{\|\delta_S\|_2 \leq r_n} \frac{1}{n}\sum_{i=1}^n \right.} - \int_0^{Z_{n,s}^{0'}\delta_S} \Big(1\{e^{0\tau}_{n,S} \leq t\} - 1\{e^{0\tau}_{n,S}\leq 0\}\Big)dt \right] \nonumber\\
	& \leq \sup_{S\in M} \mathbb{E}_{\mathcal{D}_n, \varepsilon}\left[\sup_{\|\delta_S\|_2 \leq r_n} \frac{1}{n}\sum_{i=1}^n \varepsilon_i \left(Z_{ni,S}'\delta_S\right) 1\{0\leq e^\tau_{ni,S} \leq Z_{ni,S}'\delta_S\} \right] \nonumber\\
	& \quad{} +  \sup_{S\in M} \mathbb{E}_{\mathcal{D}_n, \varepsilon}\left[\sup_{\|\delta_S\|_2 \leq r_n} \frac{1}{n}\sum_{i=1}^n \varepsilon_i e^\tau_{ni,S} 1\{0\leq e^\tau_{ni,S} \leq Z_{ni,S}'\delta_S\} \right] \nonumber\\
	&\leq \sup_{S\in M} \mathbb{E}_{\mathcal{D}_n, \varepsilon}\left[\sup_{\|\delta_S\|_2 \leq r_n} \left\|\frac{1}{n}\sum_{i=1}^n \varepsilon_i Z_{ni,S} 1\{0\leq e^\tau_{ni,S} \leq Z_{ni,S}'\delta_S\} \right\|_2 \right] r_n \label{eq:UpperBoundB1A}\\
	& \quad{} +   \sup_{S\in M} \mathbb{E}_{\mathcal{D}_n, \varepsilon}\left[\sup_{\|\delta_S\|_2 \leq r_n} \frac{1}{n}\sum_{i=1}^n \varepsilon_i e^\tau_{ni,S} 1\{0\leq e^\tau_{ni,S} \leq Z_{ni,S}'\delta_S\} \right] \label{eq:UpperBoundB1B}.
	\end{align}
	
	\textbf{Bound on eq.~\eqref{eq:UpperBoundB1A}.} Note that after de-symmetrizing eq.~\eqref{eq:UpperBoundB1A} is upper bounded by the centered quantile regression score. Thus, by the almost sure upper bound of Lemma~\ref{lemma:ASBoundQRScore},
	\begin{align*}
	&\sup_{S\in M} \mathbb{E}_{\mathcal{D}_n,\varepsilon}\left[\sup_{\|\delta_S\|_2 \leq r_n} \left\|\frac{1}{n}\sum_{i=1}^n \varepsilon_i Z_{ni,S} 1\{0\leq e^\tau_{ni,S} \leq Z_{ni,S}'\delta_S\} \right\|_2 \right] r_n \nonumber\\
	&\quad{}= O\left(\frac{1}{\lambda_n^{3/2}} \left(\frac{m \log |M| \: \log \log n}{n}\right)^{5/4} \right).
	\end{align*}
	
	\textbf{Bound on eq.~\eqref{eq:UpperBoundB1B}.} Similarly to the bound on eq~\eqref{eq:UpperBoundB1A} we conclude that
	\begin{align*}
	&\sup_{s \in M} \mathbb{E}_{\mathcal{D}_n, \varepsilon}\left[\sup_{\|\delta_S\|_2 \leq r_n} \frac{1}{n}\sum_{i=1}^n \varepsilon_i e^\tau_{ni,S} 1\{0\leq e^\tau_{ni,S} \leq Z_{ni,S}'\delta_S\} \right] \nonumber\\
	&\quad{} = O\left(\frac{1}{\lambda_n^{3/2}} \left(\frac{m\log |M| \: \log \log n}{n}\right)^{5/4} \right).
	\end{align*}
	
	\textbf{Step 3: Uniform expansion of $A_n(S)$.} 
	\begin{align}\label{eq:ExpansionA}
	&\sup_{S \in M}A_n(S) \nonumber\\
	&= \sup_{S \in M}\left|-\mathbb{E}_{\mathcal{D}_n}\left[\frac{1}{n}\sum_{i=1}^n Z_{ni,S}'\hat{\delta}_{n,S} \varphi_\tau(e^\tau_{ni,S})\right] + \mathbb{E}_{\mathcal{D}_n, (Y^0_n, X^0_n)}\left[\frac{1}{n}\sum_{i=1}^n Z_{ni}^{0'}\hat{\delta}_{n,S} \varphi_\tau(\tilde{e}^\tau_{ni,S})\right] \right. \nonumber\\
	&\quad{} \hspace{10pt}\left. + tr\left(\frac{1}{n}\sum_{i=1}^n Cov\left(Z_{ni,S}\varphi_\tau(e^\tau_{ni,S}), \: \hat{\delta}_{n,S} \right) \right)\right| \nonumber\\
	&= \sup_{S \in M}\left|-\mathbb{E}_{\mathcal{D}_n}\left[\frac{1}{n}\sum_{i=1}^n Z_{ni,S}'\hat{\delta}_{n,S} \varphi_\tau(e^\tau_{ni,S})\right] + tr \left(\frac{1}{n}\sum_{i=1}^n Cov\left(Z_{ni,S}\varphi_\tau(e^\tau_{ni,S}), \: \hat{\delta}_{n,S} \right) \right)\right| \nonumber\\
	& = 0.
	\end{align}
	
	\textbf{Step 4: Conclusion.} The claim follows by combining the upper bounds in equations~\eqref{eq:KnightsExpansion}--~\eqref{eq:ExpansionA}.
\end{proof}

\subsection{Proof of Theorem~\ref{theorem:ExpectOptTraceForm}}
\begin{proof}[\textbf{Proof of Theorem~\ref{theorem:ExpectOptTraceForm}}]
	We only need to approximate the covariance form of Theorem~\ref{theorem:ExpectOptCovForm}.
	 \begin{align*}
	 &\sup_{S \in M}\left| tr \left(\frac{1}{n}\sum_{i=1}^n Cov\left(Z_{ni,S}\varphi_\tau(e^\tau_{ni,S}), \: \hat{\delta}_{n,S} \right) \right) - \frac{1}{n}tr\left(D^\tau_{n,0}(S)^{-1}D^\tau_{n,1}(S)\right) \right| \nonumber\\
	 & = \sup_{S \in M}\left|\mathbb{E}_{\mathcal{D}_n}\left[\frac{1}{n}\sum_{i=1}^n Z_{ni,S}'\hat{\delta}_{n,S} \varphi_\tau(e^\tau_{ni,S})\right] - \frac{1}{n}tr\left(D^\tau_{n,0}(S)^{-1}D^\tau_{n,1}(S)\right) \right| \nonumber\\
	 &= \sup_{S \in M}\left|\mathbb{E}_{\mathcal{D}_n} \left[\frac{1}{n}\sum_{i=1}^n \varphi_\tau(e^\tau_{ni,S}) Z_{ni,S}' \left(D^\tau_{n,0}(S)^{-1}\frac{1}{n}\sum_{j=1}^n\varphi_\tau(e^\tau_{ni,S}) Z_{ni,S}'\right) \right] - \frac{1}{n}tr\left(D^\tau_{n,0}(S)^{-1}D^\tau_{n,1}(S)\right) \right. \nonumber\\
	 & \left. \quad{} + \mathbb{E}_{\mathcal{D}_n} \left[\frac{1}{n}\sum_{i=1}^n \varphi_\tau(e^\tau_{ni,S}) Z_{ni,S}' \left( \hat{\delta}^\tau_{ni,S} - D^\tau_{n,0}(S)^{-1} \frac{1}{n}\sum_{j=1}^n\varphi_\tau(e^\tau_{nj,S})Z_{nj,S}\right) \right] \right|\nonumber\\
	 & = \sup_{S \in M}\left|\mathbb{E}_{\mathcal{D}_n} \left[\frac{1}{n}\sum_{i=1}^n \varphi_\tau(e^\tau_{ni,S}) Z_{ni,S}' \left( \hat{\delta}^\tau_{ni,S} - D^\tau_{n,0}(S)^{-1} \frac{1}{n}\sum_{j=1}^n\varphi_\tau(e^\tau_{nj,S})Z_{nj,S}\right) \right] \right|\nonumber\\
	 & \leq \sup_{S \in M}\mathbb{E}_{\mathcal{D}_n} \left[ \left\|\frac{1}{n}\sum_{i=1}^n \varphi_\tau(e^\tau_{ni,S}) Z_{ni,S}'\right\|_2 \left\| \hat{\delta}^\tau_{ni,S} - D^\tau_{n,0}(S)^{-1} \frac{1}{n}\sum_{j=1}^n\varphi_\tau(e^\tau_{nj,S})Z_{nj,S}\right\|_2\right]\nonumber\\
	 & = O \left(c_3\left(\frac{m \log \log n}{n} \right)^{1/2} \frac{c_2}{\lambda_n^{2}}\left(\frac{m \log |M| \: \log \log n}{n}\right)^{1/4}\left(\frac{m \log |M| + \log \log n}{n}\right)^{1/2}\right) \nonumber\\
	 & = O\left(\frac{1}{\lambda_n^{2}} \left(\frac{m\log |M|\: \log \log n}{n}\right)^{5/4} \right),
	 \end{align*}
	 where the second to last equality follows from Lemmata~\ref{lemma:ASBoundQRScore} and~\ref{lemma:StrongBahadurRep}. To conclude, combine this remainder term with the one of Theorem~\ref{theorem:ExpectOptCovForm}.
\end{proof}

\subsection{Proof of Proposition~\ref{proposition:ConsistencyEstTraceForm}}
We split the proof of Proposition~\ref{proposition:ConsistencyEstTraceForm} in three parts. 

\begin{lemma}\label{lemma:ConcentrationOfD1}
	Suppose that Assumptions (A1) -- (A6) hold. Then,
	\begin{align*}
	&\sup_{S\in M}\left|tr\left(D^\tau_{n,0}(S)^{-1} \left(\widehat{D}^\tau_{n,1}(S) - D^\tau_{n,1}(S)\right)\right)  \right| = O_p\left( \frac{m}{\lambda_n^{2}} \left(\frac{m\log |M| \: \log \log n}{n}\right)^{1/2} + \frac{m}{\lambda_n}\left(\frac{\log |M|}{n}\right)^{1/2}\right).
	\end{align*}
\end{lemma}

\begin{remark}
	Since we use the quantile regression errors $\hat{e}^\tau_{ni,S}$ as proxies for the true errors $e^\tau_{ni,S}$ the process $tr\left(D^\tau_{n,0}(S)^{-1} \left(\widehat{D}^\tau_{n,1}(S) - D^\tau_{n,1}(S)\right)\right)$ is not centered. Therefore, we need to control not only the variance (standard deviation) of the process but also its deterministic drift. The deterministic drift is reflected in the first term, the variance in the second term. Note that the rate of deterministic drift can be written a $\frac{m}{\lambda_n}\times r_n$, where $r_n$ is the rate at which the estimated quantile regression vector $\hat{\theta}_{S}^\tau$converges to $\theta_{S}^\tau$ in probability, i.e. the rate at which the estimation bias of the residuals vanishes. As one expects, the rate of the term controlling the variance is proportional to the size of the maximal standard deviation (i.e. $m$) of the $n$ summands and proportional to $(\log |M|)^{1/2}$, where $|M|$ is the size of the finite set over which we take the supremum.
\end{remark}

\begin{remark}
	Clearly, under the stated assumptions, the first rate (controlling the bias) dominates the second rate (controlling the variance).
\end{remark}

\begin{proof}[\textbf{Proof of Lemma~\ref{lemma:ConcentrationOfD1}}] The goal is to apply Markov's inequality. Therefore, in the following we obtain upper bounds on the expected values of certain stochastic processes.

	\textbf{Step 1: Decomposition into deterministic bias and stochastic error terms}. By Lemmata~\ref{lemma:ASBoundQRScore} and~\ref{lemma:StrongBahadurRep} there exists $N_0 > 0$ such that for all $n \geq N_0$,
	\begin{align}\label{eq:EstimateOfD1Expansion}
	&\sup_{S \in M}\left|tr\left(D^\tau_{n,0}(S)^{-1} \left(\widehat{D}^\tau_{n,1}(S) - D^\tau_{n,1}(S)\right)\right)  \right|  \nonumber\\
	&= \sup_{S \in M} \left|\frac{1}{n} \sum_{i=1}^n \varphi_\tau^2(\hat{e}^\tau_{ni,S}) \left\| D^\tau_{n,0}(S)^{-1/2}Z_{ni,S}\right\|_2^2 - \mathbb{E}_{\mathcal{D}_n} \left[\varphi_\tau^2(e^\tau_{ni,S}) \left\| D^\tau_{n,0}(S)^{-1/2}Z_{ni,S}\right\|_2^2 \right] \right|  \nonumber\\
	& \leq \sup_{S\in M}\sup_{\|\delta_S\|_2 \leq r_n} \left|\frac{1}{n}\sum_{i=1}^n \left(\varphi_\tau^2(e^\tau_{ni,S} - Z_{ni,S}'\delta_S) - \varphi_\tau^2(e^\tau_{ni,S})\right) \left\| D^\tau_{n,0}(S)^{-1/2}Z_{ni,S}\right\|_2^2 \right. \nonumber\\
	&\phantom{\leq \sup_{S\in M}\sup_{\|\delta_S\|_2 \leq r_n} } \quad{}\left. - \mathbb{E}_{\mathcal{D}_n} \left[\left(\varphi_\tau^2(e^\tau_{ni,S} - Z_{ni,S}'\delta_S) - \varphi_\tau^2(e^\tau_{ni,S})\right)  \left\| D^\tau_{n,0}(S)^{-1/2}Z_{ni,S}\right\|_2^2 \right] \right|\nonumber\\
	& \quad{} + \sup_{S \in M}\sup_{\|\delta_S\|_2 \leq r_n} \left| \mathbb{E}_{\mathcal{D}_n} \left[\frac{1}{n} \sum_{i=1}^n \left(\varphi_\tau^2(e^\tau_{ni,S} - Z_{ni,S}'\delta_S) - \varphi_\tau^2(e^\tau_{ni,S}) \right) \left\| D^\tau_{n,0}(S)^{-1/2}Z_{ni,S}\right\|_2^2 \right] \right|\nonumber\\
	&\quad{} + \sup_{S \in M} \left|\frac{1}{n} \sum_{i=1}^n \varphi_\tau^2(e^\tau_{ni,S})\left\| D^\tau_{n,0}(S)^{-1/2}Z_{ni,S}\right\|_2^2- \mathbb{E}_{\mathcal{D}_n}\left[\varphi_\tau^2(e^\tau_{ni,S})\left\| D^\tau_{n,0}(S)^{-1/2}Z_{ni,S}\right\|_2^2\right]  \right|\nonumber\\
	& = \sup_{S\in M} A_n(S) + \sup_{S\in M} B_n(S) + \sup_{S\in M} C_n(S) \hspace{10pt} a.s.
	\end{align}
	
	\textbf{Step 2: Upper bound on $\mathbb{E}_{\mathcal{D}_n}\left[ \sup_{S\in M}  A_n(S)\right]$.} Let $\mathcal{D}^0_n$ be an independent copy of $\mathcal{D}_n$ and define
	\begin{align*}
	E_n(S) &= \mathbb{E}_{\mathcal{D}_n}\left[ \sup_{\|\delta_S\|\leq r_n}\frac{1}{n} \sum_{i=1}^n \left((1 - 2\tau)1\{0 < e^\tau_{ni,S} \leq Z_{ni,S}'\delta_S\} \left\| D^\tau_{n,0}(S)^{-1/2}Z_{ni,S}\right\|_2^2 \right. \right. \nonumber\\
	&\quad{} \hspace{90pt} \left. \left. - \mathbb{E}_{\mathcal{D}_n} \left[ (1 - 2\tau)1\{0 < e^\tau_{ni,S} \leq Z_{ni,S}'\delta_S\} \left\| D^\tau_{n,0}(S)^{-1/2}Z_{ni,S}\right\|_2^2\right] \right)\right], \nonumber\\
	W_n(S) &= \mathbb{E}_{\mathcal{D}^0_n}\left[ \sup_{\|\delta_S\|\leq r_n}\frac{1}{n^2} \sum_{i=1}^n \left((1 - 2\tau)1\{0 < e^\tau_{ni,S} \leq Z_{ni,S}'\delta\} \left\| D^\tau_{n,0}(S)^{-1/2}Z_{ni,S}\right\|_2^2 \right. \right. \nonumber\\
	&\quad{} \hspace{90pt} \left. \left. - (1 - 2\tau)1\{0 < e^{0\tau}_{ni,S} \leq Z_{ni,S}^{0'}\delta\} \left\| D^\tau_{n,0}(S)^{-1/2}Z_{ni,S}^0\right\|_2^2 \right)^2 | \mathcal{D}_n \right].
	\end{align*}
	
	Note that for $1 \leq p \leq 4$,
	\begin{align}\label{eq:BoundMoments}
	\mathbb{E}_{\mathcal{D}_n}\left[\frac{1}{n}\sum_{i=1}^n \left\|Z_{ni}\right\|_2^{2p} \right]^{1/p} \leq \mathbb{E}_{\mathcal{D}_n}\left[\frac{1}{n}\sum_{i=1}^n \left\|Z_{ni}\right\|_{2p}^{2p} \right]^{1/p} m^{1/2 - 1/(2p)} \leq \mu_{2p}^{2} m^{1/2 + 1/(2p)}.
	\end{align}
	
	By Lemma~\ref{lemma:VCMomentBound} applied to $g(Z) = (1 - 2\tau) \left\| D^\tau_{n,0}(S)^{-1/2}Z_{ni,S}\right\|_2^2$ and eq.~\eqref{eq:BoundMoments} applied to the envelope function $G(Z) = \lambda_n^{-1} \|Z\|_2^2$,
	\begin{align}\label{eq:BoundA1}
	E_n(S) \leq \frac{c_0 m}{\lambda_n}\left(\frac{m}{n}\right)^{3/4} \bigvee \frac{c_0 m}{\lambda_n} r_n^{1/2}\left(\frac{m}{n}\right)^{1/2} \leq \frac{c_0 m}{\lambda_n^{3/2}}\left(\frac{m \log |M|\:  \log \log n}{n}\right)^{3/4},
	\end{align}
	where $c_0 > 0$ is an absolute constant independent of $S \in M$.
	
	By the Hartman-Wintner law of iterated logarithm, Lemma~\ref{lemma:VCMomentBound} and eq.~\eqref{eq:BoundMoments}, there exists $N_1 > N_0$ such that for all $n \geq N_1$,
	\begin{align}\label{eq:BoundA2}
	W_n(S) &\leq \sup_{\|\delta_S\|_2\leq r_n} \frac{2}{n^2} \sum_{i=1}^n (1 - 2\tau)^2 1\{0 < e^\tau_{ni,S} \leq Z_{ni,S}'\delta\} \left\| D^\tau_{n,0}(S)^{-1/2}Z_{ni,S}\right\|_2^4 \nonumber\\
	&\quad{} + \mathbb{E}_{\mathcal{D}^0_n}\left[\sup_{\|\delta_S\|_2\leq r_n} \frac{2}{n^2} \sum_{i=1}^n (1 - 2\tau)^2 1\{0 < e^{0\tau}_{ni,S} \leq Z_{ni,S}^{0'}\delta\} \left\| D^\tau_{n,0}(S)^{-1/2}Z_{ni,S}^0\right\|_2^4\right] \nonumber\\
	&\leq \frac{2}{n^2} \sum_{i=1}^n \left(\left\| D^\tau_{n,0}(S)^{-1/2}Z_{ni,S}\right\|_2^4 - \mathbb{E}_{\mathcal{D}_n} \left[\left\| D^\tau_{n,0}(S)^{-1/2}Z_{ni,S}\right\|_2^4\right]\right)\nonumber\\
	&\quad{} +  \frac{4\nu_+}{n^2} \sum_{i=1}^n \mathbb{E}_{\mathcal{D}_n} \left[\left\| D^\tau_{n,0}(S)^{-1/2}Z_{ni,S}\right\|_2^4\right] \nonumber\\
	&\quad{} + \mathbb{E}_{\mathcal{D}^0_n}\left[\sup_{\|\delta_S\|_2\leq r_n} \frac{2}{n^2} \sum_{i=1}^n \left((1 - 2\tau)^2 1\{0 < e^{0\tau}_{ni,S} \leq Z_{ni,S}^{0'}\delta\} \left\| D^\tau_{n,0}(S)^{-1/2}Z_{ni,S}^0\right\|_2^4 \right. \right. \nonumber \\
	&\quad{} \left. \left. \hspace{110pt} - \mathbb{E}_{\mathcal{D}_n}\left[(1 - 2\tau)^2 1\{0 < e^{0\tau}_{ni,S} \leq Z_{ni,S}^{0'}\delta\} \left\| D^\tau_{n,0}(S)^{-1/2}Z_{ni,S}^0\right\|_2^4 \right] \right)\right] \nonumber\\
	&\leq \frac{c_1^2 m^2 (\log \log n)^{1/2}}{\lambda_n^2 h^2n^{3/2}} + \frac{c_1^2 m^2 r_n }{\lambda_n^2 n^{3/2}} + \frac{c_1^2 m^2 }{\lambda_n^2 n} \hspace{10pt} a.s. \nonumber\\
	&\leq \frac{c_1^2 m^2}{\lambda_n^2 n} \hspace{10pt} a.s.,
	\end{align}
	where $c_1 > 0$ is an absolute constant independent of $S \in M$.
	
	Note that
	\begin{align*}
	&\left(\varphi_\tau^2(e^\tau_{ni,S} - Z_{ni,S}'\delta) - \varphi_\tau^2(e^\tau_{ni,S})\right) \left\| D^\tau_{n,0}(S)^{-1/2}Z_{ni,S}\right\|_2^2 \\
	&\quad{}= (1 - 2\tau)1\{0 < e^\tau_{ni,S} \leq Z_{ni,S}'\delta\} \left\| D^\tau_{n,0}(S)^{-1/2}Z_{ni,S}\right\|_2^2.
	\end{align*}
	
	Thus, for fixed $S\in M$ Lemma~\ref{lemma:Panchenko} and eq.~\eqref{eq:BoundA1}--\eqref{eq:BoundA2} yield for any $t > 0$,
	\begin{align}\label{eq:BoundA3}
	&\mathbb{P}\left(\sup_{\|\delta_S\|\leq r_n}\frac{1}{n} \sum_{i=1}^n \left(\left(\varphi_\tau^2(e^\tau_{ni,S} - Z_{ni,S}'\delta) - \varphi_\tau^2(e^\tau_{ni,S})\right) \left\| D^\tau_{n,0}(S)^{-1/2}Z_{ni,S}\right\|_2^2 \right. \right.\nonumber\\
	&\quad{} \hspace{80pt} \left. \left. - \mathbb{E}_{\mathcal{D}_n} \left[\left(\varphi_\tau^2(e^\tau_{ni,S} - Z_{ni,S}'\delta) - \varphi_\tau^2(e^\tau_{ni,S})\right) \left\| D^\tau_{n,0}(S)^{-1/2}Z_{ni,S}\right\|_2^2 \right] \right) \right. \nonumber\\
	&\quad{} \hspace{10pt} \left. \geq \frac{c_0 m}{\lambda_n^{3/2}}\left(\frac{m \log |M|\:  \log \log n}{n}\right)^{3/4}  + 2  \frac{c_1 mt^{1/2}}{\lambda_n n^{1/2}} \right) \nonumber\\
	&\leq 4 e e^{-t/2}.
	\end{align}
	
	Now, for $n \geq N_1$ set $t$ to $\log|M| + t^2$ and integrate out the tail bound, 
	\begin{align}\label{eq:BoundA4}
	&\mathbb{E}_{\mathcal{D}_n}\left[ \sup_{S\in M} A_n(S)\right] \nonumber\\
	&\leq \frac{c_0 m}{\lambda_n^{3/2}}\left(\frac{m \log |M|\:  \log \log n}{n}\right)^{3/4} + \frac{2 c_1 m }{\lambda_n } \left(\frac{\log |M|}{n}\right)^{1/2} \nonumber\\
	&\quad{} + \frac{2 c_1 m }{\lambda_n n^{1/2}} \int_0^\infty \mathbb{P} \left( \sup_{S\in M} A_n(S) \geq \frac{c_0 m}{\lambda_n^{3/2}}\left(\frac{m \log |M|\:  \log \log n}{n}\right)^{3/4} +  \frac{2 c_1 m}{\lambda_n n^{1/2}} \big(t + (\log|M|)^{1/2}\big)\right) dt \nonumber\\
	&\leq \frac{c_0 m}{\lambda_n^{3/2}}\left(\frac{m \log |M|\:  \log \log n}{n}\right)^{3/4} + \frac{2 c_1 m}{\lambda_n } \left(\frac{\log |M|}{n}\right)^{1/2}  +  \frac{2 c_1 m}{\lambda_n n^{1/2}} \int_0^\infty e^{-t^2/2}dt  \nonumber\\
	&\leq \frac{c_2 m}{\lambda_n} \left(\frac{m \log |M| \: \log \log n}{n}\right)^{1/2},
	\end{align}
	where $c_2 > 0$ is an absolute constant and the last inequality follows from the rate condition (A5).

	\textbf{Step 3: Upper bound on $\sup_{S\in M} B_n(S)$.}
	\begin{align}\label{eq:BoundB1}
	\sup_{S \in M} B_n(S) &\leq 2 \sup_{S \in M} \sup_{\|\delta_S\|_2 \leq r_n} \mathbb{E}_{\mathcal{D}_n}\left[\frac{1}{n} \sum_{i=1}^n \left| F_{e^\tau_{n,S}|X_{n,S}}(Z_{ni,S}'\delta_S) - F_{e^\tau_{n,S}|X_{n,S}}(0)\right| \left\| D^\tau_{n,0}(S)^{-1/2}Z_{ni,S}\right\|_2^2 \right] \nonumber\\
	& \leq 2\nu_+ \sup_{S\in M} \sup_{\|\delta_S\|_2 \leq r_n} \mathbb{E}_{\mathcal{D}_n}\left[\frac{1}{n} \sum_{i=1}^n \left|Z_{ni,S}'\delta_S\right| \left\| D^\tau_{n,0}(S)^{-1/2}Z_{ni,S}\right\|_2^2 \right] \nonumber\\
	&\leq 2\nu_+ \sup_{S \in M} \sup_{\|\delta_S\|_2 \leq r_n} \mathbb{E}_{\mathcal{D}_n}\left[\delta_S'Z_{n,S}Z_{n,S}'\delta_S \right]^{1/2} \mathbb{E}_{\mathcal{D}_n}\left[\left\| D^\tau_{n,0}(S)^{-1/2}X_{ni,S}\right\|_2^4 \right]^{1/2} \nonumber\\
	&\leq \frac{2\lambda_+^{1/2}\nu_+}{\lambda_n} \: r_n \sup_{S \in M}  \mathbb{E}_{\mathcal{D}_n}\left[\left\|Z_{ni,S}\right\|_4^4 \right]^{1/2} m^{(4-2)/4} \nonumber\\
	&\leq \frac{2\lambda_+^{1/2}\nu_+ \mu_4^2 }{\lambda_n} \: r_n \: m \nonumber\\
	&= 2c_3 \lambda_+^{1/2}\nu_+ \mu_4^2\frac{m}{\lambda_n^2}\left(\frac{m\log|M| \: \log \log n}{n}\right)^{1/2}.	
	\end{align}
	
	\textbf{Step 4: Upper bound on $\mathbb{E}_{\mathcal{D}_n}\left[ \sup_{S\in M}  C_n(S)\right]$.} Let $\mathcal{D}^0_n$ be an independent copy of $\mathcal{D}_n$ and define
	\begin{align*}
	E_n(S) &= \mathbb{E}_{\mathcal{D}_n}\left[ \sup_{u \in \{-1,1\}} \frac{1}{n} \sum_{i=1}^n \left(\varphi_\tau^2(e^\tau_{ni,S})\left\| D^\tau_{n,0}(S)^{-1/2}Z_{ni,S}\right\|_2^2 u \right. \right.\\
	&\left. \left. \quad{} \hspace{125pt} - \mathbb{E}_{\mathcal{D}_n}\left[\varphi_\tau^2(e^\tau_{ni,S})\left\| D^\tau_{n,0}(S)^{-1/2}Z_{ni,S}\right\|_2^2 u\right] \right) \right], \nonumber\\
	W_n(S) &= \mathbb{E}_{\mathcal{D}^0_n}\left[ \sup_{u \in \{-1,1\}} \frac{1}{n^2} \sum_{i=1}^n \left(\varphi_\tau^2(e^\tau_{ni,S})\left\| D^\tau_{n,0}(S)^{-1/2}Z_{ni,S}\right\|_2^2 u \right. \right.\\
	&\left. \left. \quad{} \hspace{125pt} - \varphi_\tau^2(e^{0\tau}_{ni,S})\left\| D^\tau_{n,0}(S)^{-1/2}Z_{ni,S}^0\right\|_2^2 u  \right)^2 | \mathcal{D}_n \right].
	\end{align*}	
	Now, proceed as in Step 2. Thus, there exists $N_2 > N_1$ such that for all $n \geq N_2$,
	\begin{align}\label{eq:BoundC1}
	\mathbb{E}_{\mathcal{D}_n}\left[ \sup_{S\in M} C_n(S)\right] &\leq \frac{c_5 m}{\lambda_n} \left(\frac{\log |M|}{n}\right)^{1/2},
	\end{align}
	where $c_5 > 0$ is an absolute constant independent of $n$, $m$, and $M$.
	
	\textbf{Step 5: Conclusion:} The claim follows from Markov's inequality  and the bounds~\eqref{eq:BoundA1}--\eqref{eq:BoundC1}. Note that the bound~\eqref{eq:BoundA4} is dominated by the bound~\eqref{eq:BoundB1}. 
\end{proof}

\begin{lemma}\label{lemma:ConcentrationOfD0}
	Suppose that Assumptions (A1) -- (A6) hold. Then, for any $h > 0$,
	\begin{align*}
	&\sup_{S \in M}\left\| D^\tau_{n,0}(S)^{-1}\left( \widehat{D}_{0,h}^\tau(S) - D_{n,0}^\tau(S) \right) \right\|\\
	&\quad{}= O_p\left(\frac{h^{\alpha}}{\lambda_n} \:\bigvee \: \frac{1}{\lambda_n^{1 + \alpha}} \left(\frac{m \log |M| \:\log \log n}{n}\right)^{\alpha/2} \: \bigvee \: \frac{1}{\lambda_nh}\left(\frac{m \log |M|\: \log \log n}{n}\right)^{1/2}\right).
	\end{align*}
\end{lemma} 

\begin{remark}
	The process $\left\| D^\tau_{n,0}(S)^{-1}\left( \widehat{D}_{0,h}^\tau(S) - D_{n,0}^\tau(S) \right) \right\|$ is not centered; as in Lemma~\ref{lemma:ConcentrationOfD1} we need to control variance and a deterministic drift term: the first term captures the bias of the non-parametric estimation technique, the second term captures the bias of using the quantile regression errors $\hat{e}^\tau_{ni,S}$ as proxies for the true errors $e^\tau_{ni,S}$, and the third term captures the variance of the non-parametric estimate. Note that the rate of the second drift term can be written as $\frac{1}{\lambda_n} \times {r_n}^{\alpha}$, where $\alpha$ is the H{\"o}lder-continuity coefficient and $r_n$ is the rate at which the estimated quantile regression vector $\hat{\theta}_{S}^\tau$ converges to $\theta_{S}^\tau$ in probability, i.e. the rate at which the estimation bias of the residuals vanishes. The $\log \log n$-factor in the third term, is an artifact of our proof (for details, see comment at the beginning of the proof). However, apart from $\log \log n$-factor, the rate of the third term matches the rates of comparable results~\citep[e.g.][Theorem 5.45]{vershynin2010}.
\end{remark}

\begin{proof}[\textbf{Proof of Lemma~\ref{lemma:ConcentrationOfD0}}] The operator norm requires a different approach than the proof of Lemma~\ref{lemma:ConcentrationOfD1}. Since we take the supremum over all $S \in M$ a natural idea is to use a uniform version of Rudelson's inequality~\citep[e.g.][Lemma 3.6]{rudelson2008}. However, Rudelson's uniform inequality requires bounded predictors $Z_{ni}$ and is not easy to modify to also handle either dependent matrices or the supremum over $\delta_S \in \mathbb{R}^{|S|}$ with $\|\delta_S\|_2 \leq r_n$. Thus, instead of bounding the expected value and applying Makrov's inequality (as we did in the proof of Lemma~\ref{lemma:ConcentrationOfD1}), we use Lemma~\ref{lemma:Panchenko} to bound the tail probability, apply the union bound, and then integrate the tail probability to upper bound the expected value.
	
	Let $K_h(u) = \frac{1}{2}1\{|u| \leq h\}$.
		
	\textbf{Step 1: Decomposition into deterministic bias and stochastic error terms.} By Lemmata~\ref{lemma:ASBoundQRScore} and~\ref{lemma:StrongBahadurRep} there exists $N_0 > 0$ such that for all $n \geq N_0$,
	\begin{align*}
	&\sup_{S \in M} \left\| D^\tau_{n,0}(S)^{-1}\left( \widehat{D}_{0,h}^\tau(S) - D_{n,0}^\tau(S) \right) \right\| \nonumber\\
	&= \sup_{S \in M} \sup_{\|v\|_2=1} \left|\frac{1}{nh} \sum_{i=1}^n \left(\left|v' D^\tau_{n,0}(S)^{-1/2}Z_{ni,S}\right|^2K_h\left(e^\tau_{ni,S} - Z_{ni,S}'\hat{\delta}_{n,s}\right) \right. \right.\\
	&\left.\left. \hspace{125pt} \quad{}- h \: \mathbb{E}_{\mathcal{D}_n} \left[f_{e^\tau_{n,S}|X_{n}}(0|X_{ni}) \left|v' D^\tau_{n,0}(S)^{-1/2}Z_{ni,S}\right|^2\right]\right)\right|\\
	&\leq \sup_{S \in M} \sup_{\|v\|_2=1} \sup_{\|\delta_S\|_2 \leq r_n} \left|\frac{1}{nh} \sum_{i=1}^n \left(\left|v' D^\tau_{n,0}(S)^{-1/2}Z_{ni,S}\right|^2\left[K_h\left(e^\tau_{ni,S} - Z_{ni,S}'\delta_S\right) - K_h\left(e^\tau_{ni,S}\right) \right] \right.  \right.\\
	&\left. \left. \hspace{125pt} \quad{}  - \mathbb{E}_{\mathcal{D}_n} \left[\left|v' D^\tau_{n,0}(S)^{-1/2}Z_{ni,S}\right|^2\left[K_h\left(e^\tau_{ni,S} - Z_{ni,S}'\delta_S\right) - K_h\left(e^\tau_{ni,S}\right)  \right]\right] \right)\right| \\
	& \quad{} + \sup_{S \in M}  \sup_{\|v\|_2=1} \sup_{\|\delta_S\|_2 \leq r_n} \left|\mathbb{E}_{\mathcal{D}_n} \left[\frac{1}{nh} \sum_{i=1}^n \left|v' D^\tau_{n,0}(S)^{-1/2}Z_{ni,S}\right|^2 \left(K_h\left(e^\tau_{ni,S} - Z_{ni,S}'\delta_S\right)  - h f_{e^\tau_{n,S}|X_{n}}(0|X_{ni}) \right) \right] \right|\\
	& \quad{} + \sup_{S \in M}  \sup_{\|v\|_2=1} \left|\frac{1}{nh} \sum_{i=1}^n \left(\left|v' D^\tau_{n,0}(S)^{-1/2}Z_{ni,S}\right|^2K_h\left(e^\tau_{ni,S}\right) - \mathbb{E}_{\mathcal{D}_n}\left[\left|v' D^\tau_{n,0}(S)^{-1/2}Z_{ni,S}\right|^2K_h\left(e^\tau_{ni,S}\right)\right] \right) \right|\\
	& = \sup_{S \in M}  \sup_{\|v\|_2=1} A_n(S, v) + \sup_{S \in M} \sup_{\|v\|_2=1} B_n(S,v) + \sup_{S \in M}  \sup_{\|v\|_2=1} C_n(S,v) \hspace{10pt} a.s.
	\end{align*}
	
	\textbf{Step 2: Upper bound on $\mathbb{E}_{\mathcal{D}_n}\left[ \sup_{S \in M} \sup_{\|v\|_2=1} A_n(S, v) \right]$.} 	Let $\mathcal{D}^0_n$ be an independent copy of $\mathcal{D}_n$ and define
	\begin{align*}
	E_n(S, v) &= \mathbb{E}_{\mathcal{D}_n}\left[ \sup_{\|\delta_S\|\leq r_n}\frac{1}{nh} \sum_{i=1}^n \left(\left| v' D^\tau_{n,0}(S)^{-1/2}Z_{ni,S}\right|^2 1\{ h < e^\tau_{ni,S} \leq Z_{ni,S}'\delta_S + h \}  \right. \right. \nonumber\\
	&\quad{} \hspace{90pt} \left. \left. - \mathbb{E}_{\mathcal{D}_n} \left[ \left|v' D^\tau_{n,0}(S)^{-1/2}Z_{ni,S} \right|^2 1\{ h < e^\tau_{ni,S} \leq Z_{ni,S}'\delta_S + h \} \right] \right)\right], \nonumber\\
	W_n(S,v) &= \mathbb{E}_{\mathcal{D}^0_n}\left[ \sup_{\|\delta_S\|\leq r_n}\frac{1}{(nh)^2} \sum_{i=1}^n \left(\left| v' D^\tau_{n,0}(S)^{-1/2}Z_{ni,S}\right|^2 1\{ h < e^\tau_{ni,S} \leq Z_{ni,S}'\delta_S + h \}  \right. \right. \nonumber\\
	&\quad{} \hspace{90pt} \left. \left. - \left| v' D^\tau_{n,0}(S)^{-1/2}Z_{ni,S}^0\right|^2 1\{ h < e^{0\tau}_{ni,S} \leq Z_{ni,S}^{0'}\delta_S + h \} \right)^2 | \mathcal{D}_n \right].
	\end{align*}	
	
	By Lemma~\ref{lemma:VCMomentBound} applied to $g(Z) = \left|v'D^\tau_{n,0}(S)^{-1/2}Z_{ni,S}\right|^2$ and Lemma~\ref{lemma:MomentBounds} applied to the envelope $G(Z) \equiv g(Z)$,
	\begin{align}\label{eq:BoundAA1}
	E_n(S, v) \leq \frac{c_0}{\lambda_n h}\left(\frac{m}{n}\right)^{3/4} \bigvee \frac{c_0}{\lambda_nh } r_n^{1/2}\left(\frac{m}{n}\right)^{1/2} \leq \frac{c_0}{\lambda_n^{3/2} h}\left(\frac{m \log |M|\:  \log \log n}{n}\right)^{3/4},
	\end{align}
	where $c_0 > 0$ is an absolute constant independent of $S \in M$ and $v \in \mathbb{R}^{|S|}$.
	
	By the Hartman-Wintner law of iterated logarithm, Lemma~\ref{lemma:VCMomentBound} and Lemma~\ref{lemma:MomentBounds}, there exists $N_1 > N_0$ such that for all $n \geq N_1$,
	\begin{align}\label{eq:BoundAA2}
	W_n(S, v) &\leq \sup_{\|\delta_S\|_2\leq r_n} \frac{2}{(nh)^2} \sum_{i=1}^n \left| v' D^\tau_{n,0}(S)^{-1/2}Z_{ni,S}\right|^4 1\{ h < e^\tau_{ni,S} \leq Z_{ni,S}'\delta_S + h \} \nonumber\\
	&\quad{} + \mathbb{E}_{\mathcal{D}^0_n}\left[\sup_{\|\delta_S\|_2\leq r_n} \frac{2}{(nh)^2} \sum_{i=1}^n\left| v' D^\tau_{n,0}(S)^{-1/2}Z_{ni,S}^0\right|^4 1\{ h < e^{0\tau}_{ni,S} \leq Z_{ni,S}^{0'}\delta_S + h \}\right] \nonumber\\
	&\leq \frac{2}{(nh)^2} \sum_{i=1}^n \left(\left| v' D^\tau_{n,0}(S)^{-1/2}Z_{ni,S}\right|^4 - \mathbb{E}_{\mathcal{D}_n} \left[\left| v' D^\tau_{n,0}(S)^{-1/2}Z_{ni,S}\right|^4\right]\right)\nonumber\\
	&\quad{} + \frac{4\nu_+}{(nh)^2} \sum_{i=1}^n \mathbb{E}_{\mathcal{D}_n} \left[\left| v' D^\tau_{n,0}(S)^{-1/2}Z_{ni,S}\right|^4 \right] \nonumber\\
	&\quad{} + \mathbb{E}_{\mathcal{D}^0_n}\left[\sup_{\|\delta_S\|_2\leq r_n} \frac{2}{(nh)^2} \sum_{i=1}^n \left(\left| v' D^\tau_{n,0}(S)^{-1/2}Z_{ni,S}^0\right|^4 1\{ h < e^{0\tau}_{ni,S} \leq Z_{ni,S}^{0'}\delta_S + h \} \right. \right. \nonumber \\
	&\quad{} \left. \left. \hspace{110pt} - \mathbb{E}_{\mathcal{D}_n}\left[\left| v' D^\tau_{n,0}(S)^{-1/2}Z_{ni,S}^0\right|^4 1\{ h < e^{0\tau}_{ni,S} \leq Z_{ni,S}^{0'}\delta_S + h \}\right] \right)\right] \nonumber\\
	&\leq \frac{c_1^2(\log \log n)^{1/2}}{\lambda_n^2 h^2n^{3/2}} + \frac{c_1^2 r_n }{\lambda_n^2 h^2n^{3/2}} + \frac{c_1^2}{\lambda_n^2 h^2 n} \hspace{10pt} a.s. \nonumber\\
	&\leq \frac{c_1^2}{\lambda_n^2 h^2 n} \hspace{10pt} a.s.,
	\end{align}
	where $c_1 > 0$ is an absolute constant independent of $S \in M$ and $v \in \mathbb{R}^{|S|}$.
	
	By definition of $K_h(u)$, 
	\begin{align*}
	&\left|v' D^\tau_{n,0}(S)^{-1/2}Z_{ni,S}\right|^2\left[K_h\left(e^\tau_{ni,S} - Z_{ni,S}'\delta_S\right) - K_h\left(e^\tau_{ni,S}\right) \right]\nonumber\\
	&= \frac{1}{2}\left|v' D^\tau_{n,0}(S)^{-1/2}Z_{ni,S}\right|^2\left|1\{ h < e^\tau_{ni,S} \leq Z_{ni,S}'\delta_S + h \} - 1\{-h < e^\tau_{ni,S} \leq Z_{ni,S}'\delta_S - h\}\right|.
	\end{align*}
	
	Thus, for fixed $S\in M$ and $v \in \mathbb{R}^{|S|}$ Lemma~\ref{lemma:Panchenko} and eq.~\eqref{eq:BoundAA1}--~\eqref{eq:BoundAA2} yield for any $t > 0$,
	\begin{align}\label{eq:BoundAA3}
	&\mathbb{P}\left(\sup_{\|\delta_S\|\leq r_n}\frac{1}{nh} \sum_{i=1}^n \left(\left|v' D^\tau_{n,0}(S)^{-1/2}Z_{ni,S}\right|^2\left[K_h\left(e^\tau_{ni,S} - Z_{ni,S}'\delta_S\right) - K_h\left(e^\tau_{ni,S}\right) \right]  \right. \right.\nonumber\\
	&\quad{} \hspace{80pt} \left. \left. - \mathbb{E}_{\mathcal{D}_n} \left[\left|v' D^\tau_{n,0}(S)^{-1/2}Z_{ni,S}\right|^2\left[K_h\left(e^\tau_{ni,S}- Z_{ni,S}'\delta_S\right) - K_h\left(e^\tau_{ni,S}\right) \right] \right] \right) \right. \nonumber\\
	&\quad{} \hspace{10pt} \left. \geq \frac{c_0}{\lambda_n^{3/2} h}\left(\frac{m \log |M|\:  \log \log n}{n}\right)^{3/4}  + 2  \frac{c_1 t^{1/2}}{\lambda_n h n^{1/2}} \right) \nonumber\\
	&\leq 2 \mathbb{P}\left(\sup_{\|\delta_S\|\leq r_n}\frac{1}{nh} \sum_{i=1}^n \left(\left|v' D^\tau_{n,0}(S)^{-1/2}Z_{ni,S}\right|^2 1\{ h < e^\tau_{ni,S} \leq Z_{ni,S}'\delta_S + h \} \right. \right. \nonumber\\
	&\quad{} \hspace{95pt} \left. \left. - \mathbb{E}_{\mathcal{D}_n}\left[\left|v' D^\tau_{n,0}(S)^{-1/2}Z_{ni,S}\right|^2 1\{ h < e^\tau_{ni,S} \leq Z_{ni,S}'\delta_S + h \} \right] \right) \right. \nonumber\\
	&\quad{} \hspace{10pt} \left. \geq E_n(S, v)  + 2 W_n^{1/2}(S, v) t^{1/2}\right) \nonumber\\
	& \leq 8 e e^{-t/2}.
	\end{align}
	
	Let $\mathcal{N}_S$ be an $\frac{1}{3}$-net of the $|S|$-dimensional unit sphere. Then $|\mathcal{N}_S| \leq 7^{|S|}$ and for any symmetric $|S| \times |S|$-dimensional matrix $A$ we have $\sup_{\|v\|=1} |v'Av| \leq 3 \sup_{v \in \mathcal{N}_S}|v'Av|$~\citep[e.g.][Lemma 5.4]{vershynin2011}. Thus, $n \geq N_1$,
	\begin{align}\label{eq:BoundAA4}
	&\mathbb{E}_{\mathcal{D}_n}\left[ \sup_{S\in M} \sup_{\|v\|_2=1} A_n(S, v)\right] \nonumber\\
	&\leq 3 \mathbb{E}_{\mathcal{D}_n}\left[ \sup_{S\in M} \sup_{v \in \mathcal{N}_S} A_n(S, v)\right]\nonumber\\
	&\leq \frac{3c_0}{\lambda_n^{3/2} h}\left(\frac{m \log |M|\:  \log \log n}{n}\right)^{3/4} + \frac{6 c_1}{\lambda_n h} \left(\frac{\log |M| + \log \log n}{n}\right)^{1/2} \nonumber\\
	&\quad{} + \frac{6 c_1}{\lambda_n h n^{1/2}} \int_0^\infty \mathbb{P} \left( \sup_{S\in M} \sup_{\|v\|_2=1} A_n(S, v) \geq \frac{c_0}{\lambda_n^{3/2} h}\left(\frac{m \log |M|\:  \log \log n}{n}\right)^{3/4} \right.\nonumber\\
	&\hspace{175pt} \quad{} \left. +  \frac{2 c_1}{\lambda_n h n^{1/2}} \big(t + (\log|M| + m \log 7)^{1/2}\big)\right) dt \nonumber\\
	&\leq \frac{3c_0}{\lambda_n^{3/2} h}\left(\frac{m \log |M|\:  \log \log n}{n}\right)^{3/4} + \frac{6 c_1}{\lambda_n h} \left(\frac{\log |M| + \log \log n}{n}\right)^{1/2}  +  \frac{6 c_1}{\lambda_n hn^{1/2}} \int_0^\infty e^{-t^2/2}dt  \nonumber\\
	&\leq \frac{c_2}{\lambda_n h} \left(\frac{m \log |M| \: \log \log n}{n}\right)^{1/2},
	\end{align}
	where $c_2 > 0$ is an absolute constant and the last inequality follows from the rate condition (A5).
	
	\textbf{Step 3: Upper bound on $\sup_{S \in M}  \sup_{\|v\|_2=1} B_n(S, v)$.} The H{\"o}lder-continuity of $f_{e^\tau_{n,S}|X_{n}}$ yields
	\begin{align}\label{eq:BoundB3} 
	&\sup_{S \in M}  \sup_{\|v\|_2=1} B_n(S, v) \nonumber\\
	&\leq \sup_{S \in M}  \sup_{\|v\|_2=1} \sup_{\|\delta_S\|_2 \leq r_n} \left|\mathbb{E}_{\mathcal{D}_n} \left[\frac{1}{n} \sum_{i=1}^n \int K_1(u) \left[f_{e^\tau_{n,S}|X_{n}}(Z_{ni,S}'\delta_S + uh |X_{ni}) - f_{e^\tau_{n,S}|X_{n}}(Z_{ni,S}'\delta_S|X_{ni}) \right]du \right. \right.  \nonumber \\
	&\quad{}\left. \left. \hspace{200pt} \times \left|v' D^\tau_{n,0}(S)^{-1/2}Z_{ni,S}\right|^2\right]\right| \nonumber\\
	&\quad{}  + \sup_{S \in M} \sup_{\|v\|_2=1} \sup_{\|\delta_S\|_2 \leq r_n} \left|\mathbb{E}_{\mathcal{D}_n} \left[\frac{1}{n} \sum_{i=1}^n \int K_1(u) \left[f_{e^\tau_{n,S}|X_{n}}(Z_{ni,S}'\delta|X_{ni}) - f_{e^\tau_{n,S}|X_{n}}(0|X_{ni}) \right]du \right. \right.  \nonumber \\
	&\quad{}\left. \left. \hspace{200pt} \times  \left|v' D^\tau_{n,0}(S)^{-1/2}Z_{ni,S}\right|^2 \right] \right|\nonumber\\
	& \leq \sup_{S \in M} \sup_{\|v\|_2=1} \mathbb{E}_{\mathcal{D}_n} \left[\frac{1}{n} \sum_{i=1}^n \int K_1(u) \nu_H \left|uh\right|^{\alpha} du \left|v' D^\tau_{n,0}(S)^{-1/2}Z_{ni,S}\right|^2 \right] \nonumber\\
	& \quad{} + \sup_{S \in M} \sup_{\|v\|_2=1}  \sup_{\|\delta_S\|_2 \leq r_n} \mathbb{E}_{\mathcal{D}_n} \left[\frac{1}{n} \sum_{i=1}^n \int K_1(u) \nu_H \left|Z_{ni,S}'\delta_S\right|^{\alpha} du \left|v' D^\tau_{n,0}(S)^{-1/2}Z_{ni,S}\right|^2 \right] \nonumber\\
	&\leq \frac{\nu_H}{\lambda_n} \left(\int K_1(u) \left|u\right|^{\alpha} du\right)h^{\alpha}  + \nu_H r_n^{\alpha} \sup_{S \in M} \sup_{\|v\|_2=1}  \sup_{\|u\|_2 =1} \mathbb{E}_{\mathcal{D}_n} \left[\frac{1}{n} \sum_{i=1}^n  \left|Z_{ni,S}'u\right|^{\alpha} \left|v' D^\tau_{n,0}(S)^{-1/2}Z_{ni,S}\right|^2\right] \nonumber\\
	&\leq \frac{\nu_H}{\lambda_n} \left(\int K_1(u) \left|u\right|^{\alpha} du\right)h^{\alpha} \nonumber\\
	&\quad{} + \nu_H r_n^{\alpha} \sup_{S \in M} \sup_{\|v\|_2=1}  \sup_{\|u\|_2 =1} \left(\mathbb{E}_{\mathcal{D}_n} \left[\frac{1}{n} \sum_{i=1}^n  \left|Z_{ni,S}'u\right|^{2}\right]\right)^{\alpha/2} \left(\mathbb{E}_{\mathcal{D}_n} \left[\frac{1}{n}\sum_{i=1}^n\left|v' D^\tau_{n,0}(S)^{-1/2}Z_{ni,S}\right|^4\right]\right)^{1/2} \nonumber\\
	&\leq \frac{c_5}{\lambda_n}h^\alpha + c_\alpha r_n^{\alpha} \lambda_n^{-1}\nonumber\\
	&\leq \frac{c_5}{\lambda_n}h^\alpha + \frac{c_\alpha c_3^\alpha}{\lambda_n^{1 + \alpha}}\left(\frac{m \log |M| \: \log \log n}{n}\right)^{\alpha/2},
	\end{align}
	where $c_5, c_3, c_\alpha > 0$ are absolute constants independent of $S \in M$ and $v \in \mathbb{R}^{|S|}$ (see Lemma~\ref{lemma:MomentBounds}) .
	
	\textbf{Step 4: Upper bound on $\mathbb{E}_{\mathcal{D}_n} \left[\sup_{S \in M} \sup_{\|v\|_2=1}C_n(S, v) \right]$.} Let $\mathcal{D}^0_n$ be an independent copy of $\mathcal{D}_n$ and define
	\begin{align*}
	E_n(S, v) &= \mathbb{E}_{\mathcal{D}_n}\left[ \sup_{u \in \{-1,1\}} \frac{1}{nh} \sum_{i=1}^n \left(\left|v' D^\tau_{n,0}(S)^{-1/2}Z_{ni,S}\right|^2K_h\left(e^\tau_{ni,S}\right)u  \right.\right.\\
	&\hspace{125pt} \quad{} \left. \left. - \mathbb{E}_{\mathcal{D}_n}\left[\left|v' D^\tau_{n,0}(S)^{-1/2}Z_{ni,S}\right|^2K_h\left(e^\tau_{ni,S}\right) u\right] \right) \right], \nonumber\\
	W_n(S,v) &= \mathbb{E}_{\mathcal{D}^0_n}\left[ \sup_{u \in \{-1,1\}} \frac{1}{(nh)^2} \sum_{i=1}^n \left(\left|v' D^\tau_{n,0}(S)^{-1/2}Z_{ni,S}\right|^2K_h\left(e^\tau_{ni,S}\right) u \right.\right.\\
	&\hspace{125pt} \quad{} \left. \left. - \left|v' D^\tau_{n,0}(S)^{-1/2}Z_{ni,S}^0\right|^2K_h\left(e^{0\tau}_{ni,S}\right) u  \right)^2 | \mathcal{D}_n \right].
	\end{align*}	
	Now, proceed as in Step 2. Thus, there exists $N_2 > N_1$ such that for all $n \geq N_2$,
	\begin{align}\label{eq:BoundCC}
	\mathbb{E}_{\mathcal{D}_n}\left[ \sup_{S\in M} \sup_{\|v\|_2=1} C_n(S, v)\right] &\leq \frac{c_6}{\lambda_n h} \left(\frac{m \log |M| }{n}\right)^{1/2},
	\end{align}
	where $c_6 > 0$ is an absolute constant.
	
	\textbf{Step 5: Conclusion:} The claim follows from Markov's inequality and the bounds~\eqref{eq:BoundAA4}--~\eqref{eq:BoundCC}.
\end{proof}

\begin{proof}[\textbf{Proof of Proposition~\ref{proposition:ConsistencyEstTraceForm}}] We factor the stochastic process in two processes involving the processes in Lemmas~\ref{lemma:ConcentrationOfD1} and~\ref{lemma:ConcentrationOfD0}. Then, convergence in probability at the given rate follows immediately.
	
	\textbf{Factorization.}
	\begin{align*}
	&\sup_{S \in M} \left|tr\left(\widehat{D}_{0,h}^\tau(S)^{-1}\widehat{D}_{n,1}^\tau(S)\right) - tr\left(D^\tau_{n,0}(S)^{-1}D^\tau_{n,1}(S)\right) \right|\\
	& = \sup_{S \in M} \left| tr \left(\left(\widehat{D}_{0,h}^\tau(S)^{-1} - D^\tau_{n,0}(S)^{-1}\right)\widehat{D}_{n,1}^\tau(S) \right) + tr\left(D^\tau_{n,0}(S)^{-1}\left(\widehat{D}_{n,1}^\tau(S) - D^\tau_{n,1}(S)\right)\right)\right|\\
	&\leq \left(\sup_{S \in M} \left\|\left(D_{0,S} - \widehat{D}_{0,S}^h\right) D^\tau_{n,0}(S)^{-1}\right\| \right) \left(\sup_{S \in M} \left|tr\left(\widehat{D}_{n,1}^\tau(S) \widehat{D}_{0,h}^\tau(S)^{-1}\right)\right|\right) \\
	& \quad{} + \sup_{S \in M} \left|tr\left(D^\tau_{n,0}(S)^{-1}\left(\widehat{D}_{n,1}^\tau(S) - D^\tau_{n,1}(S)\right)\right)\right| \\
	&\leq \left(\sup_{S \in M} \left\|\left(D_{0,S} - \widehat{D}_{0,S}^h\right) D^\tau_{n,0}(S)^{-1}\right\| \right) \left(\sup_{S \in M} \left|tr\left(D_{n,1}^\tau(S) D_{0,h}^\tau(S)^{-1}\right)\right|\right)\\
	&\quad{} + \left(\sup_{S \in M} \left\|\left(D_{0,S} - \widehat{D}_{0,S}^h\right) D^\tau_{n,0}(S)^{-1}\right\| \right) \left(\sup_{S \in M} \left|tr\left( \widehat{D}_{n,1}^\tau(S) \widehat{D}_{0,h}^\tau(S)^{-1} - D_{n,1}^\tau(S) D_{0,h}^\tau(S)^{-1}\right)\right|\right)\\
	&\quad{} + \sup_{S \in M} \left|tr\left(D^\tau_{n,0}(S)^{-1}\left(\widehat{D}_{n,1}^\tau(S) - D^\tau_{n,1}(S)\right)\right)\right|.
	\end{align*}
	Re-arranging and solving for $\sup_{S \in M} \left|tr\left(\widehat{D}_{0,h}^\tau(S)^{-1}\widehat{D}_{n,1}^\tau(S)\right) - tr\left(D^\tau_{n,0}(S)^{-1}D^\tau_{n,1}(S)\right) \right|$ gives
	\begin{align}
	\begin{split}
	&\sup_{S \in M} \left|tr\left(\widehat{D}_{0,h}^\tau(S)^{-1}\widehat{D}_{n,1}^\tau(S)\right) - tr\left(D^\tau_{n,0}(S)^{-1}D^\tau_{n,1}(S)\right) \right| \\
	& = \frac{\sup_{S \in M} \left\|\left(D_{0,S} - \widehat{D}_{0,S}^h\right) D^\tau_{n,0}(S)^{-1}\right\| }{1- \sup_{S \in M} \left\|\left(D_{0,S} - \widehat{D}_{0,S}^h\right) D^\tau_{n,0}(S)^{-1}\right\| } \sup_{S \in M} \left|tr\left(D_{n,1}^\tau(S) D_{0,h}^\tau(S)^{-1}\right)\right| \\
	&\quad{} + \frac{\sup_{S \in M} \left|tr\left(D^\tau_{n,0}(S)^{-1}\left(\widehat{D}_{n,1}^\tau(S) - D^\tau_{n,1}(S)\right)\right)\right|}{1- \sup_{S \in M} \left\|\left(D_{0,S} - \widehat{D}_{0,S}^h\right) D^\tau_{n,0}(S)^{-1}\right\| }.
	\end{split}
	\end{align}
	Thus, by Lemmata~\ref{lemma:ConcentrationOfD1} and~\ref{lemma:ConcentrationOfD0},
	\begin{align}
	&\sup_{S \in M} \left|tr\left(\widehat{D}_{0,h}^\tau(S)^{-1}\widehat{D}_{n,1}^\tau(S)\right) - tr\left(D^\tau_{n,0}(S)^{-1}D^\tau_{n,1}(S)\right) \right| \nonumber\\
	&= O_p\left(\frac{m h^{\alpha}}{\lambda_n^2} \:\bigvee \: \frac{m}{\lambda_n^{2 + \alpha}} \left(\frac{m \log |M| \:\log \log n}{n}\right)^{\alpha/2} \: \bigvee \: \frac{m}{\lambda_n^2h}\left(\frac{m \log |M|\: \log \log n}{n}\right)^{1/2}\right) \nonumber\\
	& + O_p\left( \frac{m}{\lambda_n^{2}} \left(\frac{m\log |M| \: \log \log n}{n}\right)^{1/2} + \frac{m}{\lambda_n}\left(\frac{\log |M|}{n}\right)^{1/2}\right) \nonumber\\
	&= O_p\left(\frac{m h^{\alpha}}{\lambda_n^2} + \frac{m}{\lambda_n^2h}\left(\frac{m \log |M|\: \log \log n}{n}\right)^{1/2} + \frac{m}{\lambda_n^{2 + \alpha}} \left(\frac{m \log |M| \:\log \log n}{n}\right)^{\alpha/2}\right),
	\end{align} 
	by Assumption (A5) on the lower bound on $\lambda_n$.

\end{proof}

\subsection{Proof of Theorem~\ref{theorem:ConsistencyEstExpectOpt}}
\begin{proof}[\textbf{Proof of Theorem~\ref{theorem:ConsistencyEstExpectOpt}}]
	
	By Theorem~\ref{theorem:ExpectOptTraceForm} we have $\inf_{S \in M} b_n^\tau(S) > 0$ and
	\begin{align}\label{eq:LowerBoundOpt}
	\inf_{S \in M} b_n^\tau(S) &\gtrsim O\left(n^{-1}\right).
	\end{align}
	Let $r_{n,2}$ be as defined in Theorem~\ref{theorem:ExpectOptTraceForm} and fix $T > 0$. By Theorem~\ref{theorem:ExpectOptTraceForm} there exists $T > 0$  and $N > 0$ such that for all $n \geq N$, 
	\begin{align*}
	\sup_{S \in M} \left| b_n^\tau(S)- \frac{1}{n}tr\left(D^\tau_{n,0}(S)^{-1}D^\tau_{n,1}(S)\right)\right| \leq T r_{n,2}.
	\end{align*}
	Thus, by Proposition~\ref{proposition:ConsistencyEstTraceForm} for any $\epsilon > 0$ there exist $T' > T$ and $N' \geq N$ such that for all $n > N'$,
	\begin{align*}
	&\mathbb{P}\left( \sup_{S \in M} \left|\frac{\hat{b}_n^\tau(S)}{b_n^\tau(S)} - 1 \right| > \frac{T}{\inf_{S \in M} b_n^\tau(S)} \left(\frac{m\: h^\alpha}{\lambda_n n} + \frac{m\: r_{n}}{n h} + \frac{m \: r_{n}^\alpha}{\lambda_n n} + r_{n,2} \right) \right) \nonumber\\
	& \leq \mathbb{P}\left(\sup_{S \in M} \left|\hat{b}_n^\tau(S)- \frac{1}{n}tr\left(D^\tau_{n,0}(S)^{-1}D^\tau_{n,1}(S)\right)\right| + \sup_{S \in M} \left| b_n^\tau(S)- \frac{1}{n}tr\left(D^\tau_{n,0}(S)^{-1}D^\tau_{n,1}(S)\right)\right| \right. \nonumber\\
	&\quad{} \left. \hspace{50pt}> \frac{T}{n} \left(\frac{h^\alpha}{\lambda_n} + \frac{r_{n}}{h} + \frac{m\: r_{n}^\alpha}{\lambda_n} \right) + T r_{n,2} \right) \nonumber\\
	&\leq \mathbb{P}\left(\sup_{S \in M} \left|n\hat{b}_n^\tau(S)- tr\left(D^\tau_{n,0}(S)^{-1}D^\tau_{n,1}(S)\right)\right| >  T\left(\frac{m\: h^\alpha}{\lambda_n} + \frac{m\: r_{n}}{h} + \frac{m \: r_{n}^\alpha}{\lambda_n}\right) \right) \nonumber\\
	& < \epsilon.
	\end{align*}
	
	To conclude, note that by eq.~\eqref{eq:LowerBoundOpt}
	\begin{align*}
	\frac{r_{n,4}}{\inf_{S \in M} b_n^\tau(S)}  &= O\left( \frac{n}{\lambda_n} \left(\frac{m \log|M| \: \log \log n}{n}\right)^{5/4} + \left(\frac{m\: h^\alpha}{\lambda_n} + \frac{m\: r_{n}}{h} + \frac{m \: r_{n}^\alpha}{\lambda_n}\right) \right)\\
	&= O\left( n \lambda_n^{3/2} r_{n}^{5/2} +  \frac{m\: h^\alpha}{\lambda_n} + \frac{m\: r_{n}}{h} + \frac{m \: r_{n}^\alpha}{\lambda_n}\right).
	\end{align*}
\end{proof}

\subsection{Proof of Theorem~\ref{theorem:PredRiskConsistency}}
\begin{proof}[\textbf{Proof of Theorem~\ref{theorem:PredRiskConsistency}}]
	\textbf{Step 1: Decomposition.} 
	\begin{align}
	&\sup_{S \in M}\left|\frac{1}{n}\sum_{i=1}^n\rho_\tau(Y_{ni} - Z_{ni,S}'\hat{\theta}^\tau_{n,S}) - \rho_\tau(Y_{ni}) + tr\left(\widehat{D}^\tau_{h,0}(S)^{-1}\widehat{D}^\tau_{n,1}(S)\right) \right. \nonumber\\
	& \quad{} \hspace{30pt} \left. - \mathbb{E}_{\mathcal{D}_n, (Y^0_n, X^0_n)} \left[\frac{1}{n}\sum_{i=1}^n \rho_\tau(Y^0_{ni} - Z_{ni,S}^{0'}\hat{\theta}^\tau_{n,S}) - \rho_\tau(Y^0_{ni})\right] \right| \nonumber\\
	& \leq \sup_{S \in M}\left| \frac{1}{n}\sum_{i=1}^n\rho_\tau(Y_{ni} - Z_{ni,S}'\hat{\theta}^\tau_{n,S}) - \rho_\tau(Y_{ni} - Z_{ni,S}'\theta^\tau_{n,S}) \right. \nonumber\\
	&\quad{} \hspace{30pt} \left. - \mathbb{E}_{\mathcal{D}_n} \left[\frac{1}{n}\sum_{i=1}^n\rho_\tau(Y_{ni} - Z_{ni,S}'\hat{\theta}^\tau_{n,S}) - \rho_\tau(Y_{ni} - Z_{ni,S}'\theta^\tau_{n,S})\right] \right| \nonumber\\
	& \quad{} + \sup_{S \in M}\left| \mathbb{E}_{\mathcal{D}_n} \left[\frac{1}{n}\sum_{i=1}^n\rho_\tau(Y_{ni} - Z_{ni,S}'\hat{\theta}^\tau_{n,S}) - \rho_\tau(Y_{ni})\right] - \mathbb{E}_{\mathcal{D}_n, (Y^0_n, X^0_n)} \left[\frac{1}{n}\sum_{i=1}^n \rho_\tau(Y^0_{ni} - Z_{ni,S}^{0'}\hat{\theta}^\tau_{n,S}) - \rho_\tau(Y^0_{ni})\right] \right. \nonumber\\
	&\quad{} \hspace{30pt} \left. + \frac{1}{n}tr\left(D^\tau_{n,0}(S)^{-1}D^\tau_{n,1}(S)\right)  \right| \nonumber\\
	& \quad{} + \sup_{S \in M}\frac{1}{n}\left|  tr\left(\widehat{D}^\tau_{h,0}(S)^{-1}\widehat{D}^\tau_{n,1}(S)\right) - tr\left(D^\tau_{n,0}(S)^{-1}D^\tau_{n,1}(S)\right) \right| \nonumber\\
	&\quad{} + \sup_{S \in M}\left| \frac{1}{n}\sum_{i=1}^n\rho_\tau(Y_{ni} - Z_{ni,S}'\theta^\tau_{n,S}) - \rho_\tau(Y_{ni})  -  \mathbb{E}_{\mathcal{D}_n} \left[\frac{1}{n}\sum_{i=1}^n\rho_\tau(Y_{ni} - Z_{ni,S}'\theta^\tau_{n,S}) - \rho_\tau(Y_{ni})\right] \right| \nonumber\\
	& = \sup_{S \in M}A_n(S) + \sup_{S \in M}B_n(S) + \sup_{S \in M}C_n(S) + \sup_{S \in M}D_n(S).
	\end{align}
	
	\textbf{Step 2: Bounds on $\sup_{S \in M}B_n(S)$ and $\sup_{S \in M}C_n(S)$.} By Theorem~\ref{theorem:ExpectOptTraceForm},
	\begin{align*}
	\sup_{S \in M}B_n(S)= O \left(\lambda_n^{3/2} r_{n}^{5/2}\right) \hspace{10pt} a.s.,
	\end{align*}
	and by Theorem~\ref{theorem:ConsistencyEstExpectOpt}, 
	\begin{align*}
	\sup_{S \in M}C_n(S) &= O_p \left( \frac{m \: h^\alpha}{\lambda_n^2 n} + \frac{m\:r_{n}}{h \lambda_n n} + \frac{m\: r_{n}^\alpha}{\lambda_n^2 n} \right).
	\end{align*}
	
	\textbf{Step 3: Bound on $\sup_{S \in M}D_n(S)$.} Recall that $\theta^\tau_{n,S}$ solves the population quantile regression minimization problem under the constraint that the minimizer is a linear function, while $Q_{Y_n}(\tau|X_n)$ solves the unconstrained population quantile regression minimization problem. Therefore, we have for all $i =1, \ldots, n$,
	\begin{align}\label{eq:MeanBound}
	0 &\leq \mathbb{E}_{\mathcal{D}_n}\left[ \min_{S \in M}\rho_\tau(Y_{ni}) - \rho_\tau(Y_{ni} - Z_{ni,S}'\theta^\tau_{n,S} \right] \nonumber\\
	& \leq \mathbb{E}_{\mathcal{D}_n}\left[\min_{S \in M} \rho_\tau(Y_{ni}) - \rho_\tau(Y_{ni} - Q_{Y_n}(\tau|X_{ni})) \right] \nonumber\\
	& \leq \mathbb{E}_{\mathcal{D}_n}\left[\big| Q_{Y_n}(\tau|X_{ni}))\big| \right] < \infty.
	\end{align}
	
	The chain of inequalities in~\ref{eq:MeanBound}, the convexity of the maximum operator together with Jensen's Inequality, and eq.~\ref{eq:PopulationQRProblemWeightedSquare} imply
	\begin{align*}
	0 &\leq \mathbb{E}_{\mathcal{D}_n}\left[\min_{S\in M} \left\{\rho_\tau(Y_{ni}) - \rho_\tau(Y_{ni} - Z_{ni,S}'\theta^\tau_{n,S})\right\} \right]\\
	&\leq \mathbb{E}_{\mathcal{D}_n}\left[\big| Q_{Y_n}(\tau|X_{ni}))\big| \right]  + \mathbb{E}_{\mathcal{D}_n}\left[\min_{S \in M} \left\{ - \rho_\tau(Y_{ni}-Z_{ni,S}'\theta^\tau_{n,S} ) + \rho_\tau(Y_{ni} - Q_{Y_n}(\tau|X_{ni}))\right\} \right]\\
	&\leq \mathbb{E}_{\mathcal{D}_n}\left[\big| Q_{Y_n}(\tau|X_{ni}))\big| \right]  - \mathbb{E}_{\mathcal{D}_n}\left[\max_{S \in M} \left\{ \rho_\tau(Y_{ni}-Z_{ni,S}'\theta^\tau_{n,S} ) - \rho_\tau(Y_{ni} - Q_{Y_n}(\tau|X_{ni}))\right\} \right]\\
	& \leq \mathbb{E}_{\mathcal{D}_n}\left[\big| Q_{Y_n}(\tau|X_{ni}))\big| \right]  - \frac{\nu_-}{2}\mathbb{E}_{\mathcal{D}_n}\left[\max_{S \in M} \left( Z_{ni,S}'\theta^\tau_{n,S} - Q_{Y_n}(\tau|X_n)\right)^2\right].
	\end{align*}
	For $i = 1, \ldots, n$ above inequality gives
	\begin{align}\label{eq:UnifUpperBoundExpApproxError}
	\mathbb{E}_{\mathcal{D}_n}\left[\max_{S \in M} \left( Z_{ni,S}'\theta^\tau_{n,S} - Q_{Y_n}(\tau|X_{ni})\right)^2\right] \leq \frac{2}{\nu_-}\mathbb{E}_{\mathcal{D}_n}\left[\big| Q_{Y_n}(\tau|X_{ni}))\big| \right] < \infty.
	\end{align}
	
	Let $\varepsilon=(\varepsilon_1, \ldots, \varepsilon_n)$ be a vector of i.i.d. Rademacher variables independent of $\mathcal{D}_n$. By the Sub-Gaussianity of the conditional Rademacher average, the Lipschitz continuity of the quantile loss function, and eq.~\eqref{eq:UnifUpperBoundExpApproxError},
	\begin{align*}
	&\mathbb{E}_{\mathcal{D}_n}\left[\sup_{S \in M}\left|\frac{1}{n}\sum_{i=1}^n \rho_\tau(Y_{ni} - Z_{ni,S}\theta^\tau_{n,S}) - \rho_\tau(Y_{ni}) - \mathbb{E}_{\mathcal{D}_n}\left[\rho_\tau(Y_{ni} - Z_{ni,S}\theta^\tau_{n,S}) - \rho_\tau(Y_{ni})\right] \right| \right]\\
	&\leq \left(\mathbb{E}_{\mathcal{D}_n} \left[\sup_{S \in M}\frac{1}{n}\sum_{i=1}^n \left(Z_{ni,S}'\theta^\tau_{n,S}\right)^2 \right]\right)^{1/2}\left(\frac{\log |M|}{n}\right)^{1/2}\\
	& \leq \left[\left(\mathbb{E}_{\mathcal{D}_n} \left[\sup_{S \in M}\frac{1}{n}\sum_{i=1}^n \left(Z_{ni,S}'\theta^\tau_{n,S} - Q_{Y_n}(\tau|X_{ni})\right)^2 \right]\right)^{1/2} + \left(\mathbb{E}_{\mathcal{D}_n} \left[\frac{1}{n}\sum_{i=1}^n Q^2_{Y_n}(\tau|X_{ni}) \right]\right)^{1/2}\right]\left(\frac{\log |M|}{n}\right)^{1/2}\\
	&\leq c_1\left(\frac{\log |M|}{n}\right)^{1/2},
	\end{align*}
	where $c_1 >0$ depends on the constants in eq.~\eqref{eq:UnifUpperBoundExpApproxError}. Thus,
	\begin{align*}
	\sup_{S \in M}D_n(S) = O_p\left(\left(\frac{\log |M|}{n}\right)^{1/2}\right).
	\end{align*}
	
	\textbf{Step 4: Bound on $\sup_{S \in M}A_n(S)$}. 
	Note that 
	\begin{align}\label{eq:DecompositionRho}
	\begin{split}
	2\Big(\rho_\tau(Y- Z'\theta_1) - \rho_\tau(Y- Z'\theta_2) \Big)&= Z'(\theta_1-\theta_2) 1\{Y \geq Z'\theta_1 \}1\{Y \geq Z'\theta_2\}\\
	&\quad{} - Z'(\theta_1-\theta_2) 1\{Y < Z'\theta_1 \}1\{Y < Z'\theta_2\} \\
	&\quad{} + (2Y- Z'\theta_1-Z'\theta_2) 1\{Y \geq Z'\theta_1 \}1\{Y < Z'\theta_2\}\\
	&\quad{} - (2Y- Z'\theta_1-Z'\theta_2) 1\{Y < Z'\theta_1 \}1\{Y \geq Z'\theta_2\} \\
	&\quad{} + (2\tau-1)Z'(\theta_1 - \theta_2).
	\end{split}
	\end{align}
	Let $\mathcal{D}^0_n$ be an independent copy of $\mathcal{D}_n$ and define
	\begin{align*}
	W_n &= \mathrm{E}_{\mathcal{D}^0_n}\left[\max_{S \in M} \sup_{\|\delta_S\|_2 \leq r_n}\frac{1}{n^2} \sum_{i=1}^n \Big(\rho_\tau(e^\tau_{ni,S}- Z_{ni,S}'\delta_S) - \rho_\tau(e^\tau_{ni,S})  - \rho_\tau(e^{0\tau}_{ni,S}- Z_{ni,S}^{0'}\delta_S) - \rho_\tau(e^{0\tau}_{ni,S}) \Big)^2 | \mathcal{D}_n \right].
	\end{align*}
	
	By expansion~\eqref{eq:DecompositionRho},
	\begin{align}\label{eq:WnConsistency}
	W_n &\leq \max_{S \in M} \sup_{\|\delta_S\|_2 \leq r_n}\frac{2}{n^2} \sum_{i=1}^n \Big(\rho_\tau(e^\tau_{ni,S}- Z_{ni,S}'\delta_S) - \rho_\tau(e^\tau_{ni,S})\Big)^2 \nonumber\\
	&\quad{} + \mathrm{E}\left[\max_{S \in M} \sup_{\|\delta_S\|_2 \leq r_n}\frac{2}{n^2}\sum_{i=1}^n  \Big(\rho_\tau(e^{0\tau}_{ni,S}- Z_{ni,S}^{0'}\delta_S) - \rho_\tau(e^{0\tau}_{ni,S}) \Big)^2\right] \nonumber\\
	&\leq \max_{S \in M} \sup_{\|\delta_S\|_2 \leq r_n}\frac{9}{n^2}\sum_{i=1}^n \left(\big(Z_{ni,S}'\delta_S\big)^2  - \mathbb{E}_{\mathcal{D}_n}\left[\big(Z_{ni,S}'\delta_S\big)^2\right] \right) \nonumber\\
	&\quad{} + \mathbb{E}_{\mathcal{D}^0_n}\left[\max_{S \in M} \sup_{\|\delta_S\|_2 \leq r_n}\frac{9}{n^2}\sum_{i=1}^n \left(\big(Z_{ni,S}^{0'}\delta_S\big)^2  - \mathbb{E}_{\mathcal{D}_n}\left[\big(Z_{ni,S}^{0'}\delta_S\big)^2\right] \right) \right] \nonumber\\
	&\quad{} + \max_{S \in M} \sup_{\|\delta_S\|_2 \leq r_n} \mathbb{E}_{\mathcal{D}_n}\left[\frac{18}{n^2}\sum_{i=1}^n \big(Z_{ni,S}'\delta_S\big)^2\right] \hspace{10pt} a.s.
	\end{align}
	
	As in Step 2 of the proof of Lemma~\ref{lemma:ASBoundQRScore},  we conclude that there exist absolute constants $N_2, c_2 > 0$ such that for all $n \geq N_2$,
	\begin{align*}
	\max_{S \in M} \sup_{\|\delta_S\|_2 \leq r_n}\frac{9}{n^2}\sum_{i=1}^n \left(\big(Z_{ni,S}'\delta_S\big)^2  - \mathbb{E}_{\mathcal{D}_n}\left[\big(Z_{ni,S}'\delta_S\big)^2\right] \right) \leq c_2 \left(\frac{m \log|M| \log \log n}{n}\right)^{1/2}\frac{r_n^2}{n} \hspace{10pt} a.s.
	\end{align*}
	Thus, there exist absolute constants $N_3, c_3 > 0$ such that for all $n \geq N_3$,
	\begin{align}\label{eq:BoundAAA1}
	W_n \leq c_3 \frac{r_n^2}{n} \hspace{10pt} a.s.
	\end{align}
	
	By the second statement of Lemma~\ref{lemma:Panchenko} for any $t > 0$ and all $n \geq N_3$, 
	\begin{align}\label{eq:BoundAAA2}
	&\mathbb{P}\left(\sup_{S \in M} \frac{1}{n}\sum_{i=1}^n \Big(\rho_\tau(Y_{ni} - Z_{ni,S}'\hat{\theta}^\tau_{n,S}) - \rho_\tau(Y_{ni} - Z_{ni,S}'\theta^\tau_{n,S}) \Big) \right. \nonumber\\
	&\quad{} \hspace{30pt} \left. - \mathbb{E}_{\mathcal{D}_n} \left[\frac{1}{n}\sum_{i=1}^n\rho_\tau(Y_{ni} - Z_{ni,S}'\hat{\theta}^\tau_{n,S}) - \rho_\tau(Y_{ni} - Z_{ni,S}'\theta^\tau_{n,S})\right] \geq 2 W_n^{1/2}t\right) \nonumber\\
	&\leq \sum_{S \in M} \mathbb{P}\left(- \sup_{\|\delta_S\|_2 \leq r_n} -\frac{1}{n}\sum_{i=1}^n \Big(\rho_\tau(Y_{ni} - Z_{ni,S}'\hat{\theta}^\tau_{n,S}) - \rho_\tau(Y_{ni} - Z_{ni,S}'\theta^\tau_{n,S}) \Big) \right. \nonumber\\
	&\quad{} \hspace{30pt} \left. + \: \mathbb{E}_{\mathcal{D}_n} \left[\sup_{\|\delta_S\|_2 \leq r_n} - \frac{1}{n}\sum_{i=1}^n\rho_\tau(Y_{ni} - Z_{ni,S}'\hat{\theta}^\tau_{n,S}) - \rho_\tau(Y_{ni} - Z_{ni,S}'\theta^\tau_{n,S})\right] \geq 2 W_n^{1/2}t\right) \nonumber\\
	&\leq \sum_{S \in M} \mathbb{P}\left( \sup_{\|\delta_S\|_2 \leq r_n} -\frac{1}{n}\sum_{i=1}^n \Big(\rho_\tau(Y_{ni} - Z_{ni,S}'\hat{\theta}^\tau_{n,S}) - \rho_\tau(Y_{ni} - Z_{ni,S}'\theta^\tau_{n,S})\Big) \right. \nonumber\\
	&\quad{} \hspace{30pt} \left. \leq \mathbb{E}_{\mathcal{D}_n} \left[\sup_{\|\delta_S\|_2 \leq r_n} - \frac{1}{n}\sum_{i=1}^n\rho_\tau(Y_{ni} - Z_{ni,S}'\hat{\theta}^\tau_{n,S}) - \rho_\tau(Y_{ni} - Z_{ni,S}'\theta^\tau_{n,S})\right] - 2 W_n^{1/2}t\right) \nonumber\\
	&\leq 4 |M| e e^{-t^2/2}.
	\end{align}
	Analogously, we derive a bound on the probability for the lower tail via the first statement of Lemma~\ref{lemma:Panchenko},
	\begin{align}\label{eq:BoundAAA3}
	&\mathbb{P}\left(\sup_{S \in M} \frac{1}{n}\sum_{i=1}^n \Big(\rho_\tau(Y_{ni} - Z_{ni,S}'\hat{\theta}^\tau_{n,S}) - \rho_\tau(Y_{ni} - Z_{ni,S}'\theta^\tau_{n,S}) \Big) \right. \nonumber\\
	&\quad{} \hspace{30pt} \left. - \mathbb{E}_{\mathcal{D}_n} \left[\frac{1}{n}\sum_{i=1}^n\rho_\tau(Y_{ni} - Z_{ni,S}'\hat{\theta}^\tau_{n,S}) - \rho_\tau(Y_{ni} - Z_{ni,S}'\theta^\tau_{n,S})\right] \leq -2 W_n^{1/2}t\right) \nonumber\\
	&\leq 4 |M| e e^{-t^2/2}.
	\end{align}
	Thus, combining eq.~\eqref{eq:BoundAAA1}-~\eqref{eq:BoundAAA3} and setting $t = t'(\log|M|)^{1/2}$, there exists $N_4 > N_3$ such that for all $n > N_4$,
	\begin{align*}
	\mathrm{P}\left( \sup_{S \in M} A_n(S) \geq c_4 r_n\frac{t^{1/2}}{n^{1/2}} \right) \leq 4 e e^{-t'/2},
	\end{align*}
	where $c_4 > 0$ is an absolute constant. 
	Hence, as in Step 2 of the proof of Lemma~\ref{lemma:ConcentrationOfD0} integrating this tail bound out yields for all $ n > N_4$, $\mathbb{E}_{\mathcal{D}_n}\left[\sup_{S \in M} A_n(S)\right] \leq c_5 \frac{r_n}{n^{1/2}}$, where $c_5 > 0$ is an absolute constant, and
	\begin{align*}
	\mathbb{E}_{\mathcal{D}_n}\left[\sup_{S \in M} A_n(S)\right] = O_p\left(\frac{r_n}{n^{1/2}}\right).
	\end{align*}
	
	\textbf{Step 5: Conclusion.} Combining above bounds on $\sup_{S \in M}A_n(S)$, $\sup_{S \in M}B_n(S)$, $\sup_{S\in M}C_n(S)$, and $\sup_{S \in M}D_n(S)$  we have
	\begin{align*}
	\sup_{S \in M}\left| \widehat{PR}^{\tau}_{n,h}(S) - \mathrm{PR}^\tau_n(S) \right| =O_p\left(\left(\frac{\log |M|}{n}\right)^{1/2} + \frac{r_n}{n^{1/2}} + \lambda_n^{3/2} r_{n,3}^{5/2} + \frac{m \: h^\alpha}{\lambda_n^2 n} + \frac{m\:r_{n,3}}{h \lambda_n n} + \frac{m\: r_{n,3}^\alpha}{\lambda_n^2 n} \right).
	\end{align*}	
\end{proof}

\includepdf[pages={{},-}]{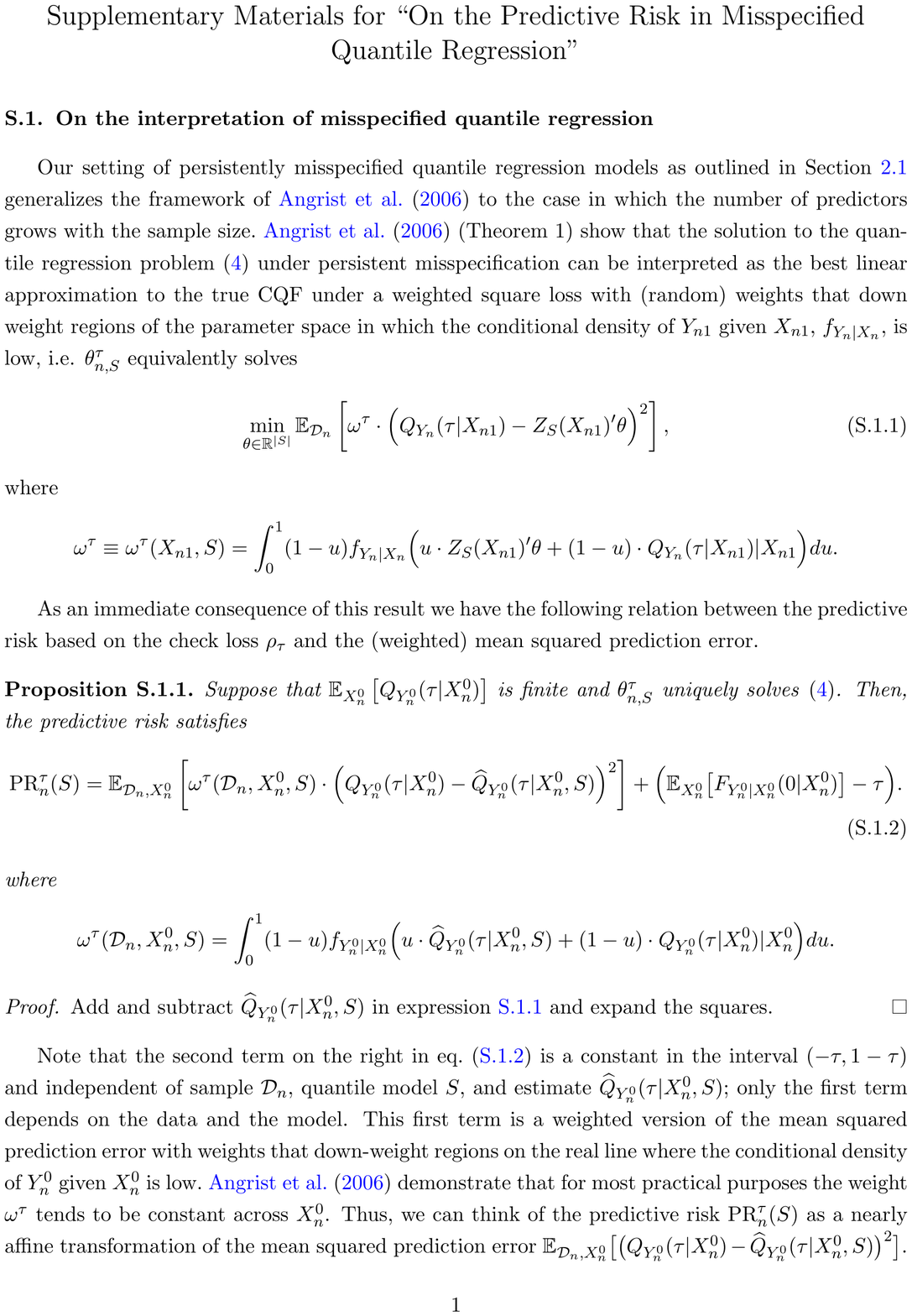}

\end{document}